\documentclass[12pt,reqno,english,a4]{amsart}
\usepackage{amsmath,amsthm,amsfonts,amssymb,amscd}
\usepackage[latin1]{inputenc}
\usepackage{psfrag}
\usepackage{epsfig}
\usepackage{a4wide}
\usepackage[]{graphicx}

\usepackage{hyperref} % criar link
\theoremstyle{plain}
\newtheorem{theorem}{Theorem}

\newtheorem{proposition}{Proposition}
\newtheorem{lemma}{Lemma}

\newtheorem{example}{Example}
\newtheorem{remark}{Remark}
\newtheorem{corollary}{Corollary}
\newtheorem{definition}{Definition}

\numberwithin{equation}{section}
\renewcommand{\th} {\theta}

\newcommand{\vphi}{\varphi}
\newcommand{\om} {\omega}       
\newcommand{\fin}{\hfill$\square$}

\def \CC {{\mathbb C}}
\def \RR {{\mathbb R}}
\def \ZZ {{\mathbb Z}}
\def \NN {{\mathbb N}}

\def \TT {{\mathbb T}}
\def \QQ {{\mathbb Q}}
\def \SS {{\mathbb S}}
\def \cO {{\mathcal O}}
\def \cF {{\mathcal F}}
\def \cG {{\mathcal G}}

\def \cR {{\mathcal R}}

\def\cA{{\mathcal A}}
\def\cB{{\mathcal B}}
\def\cC{{\mathcal C}}
\def\cH{{\mathcal H}}
\def\cS{{\mathcal S}}

\def\-{{\setminus}}

\def\vphi{\varphi}

\begin{document}

\title[Transitivity of codimension one Anosov actions of $\RR^k$]
      {Transitivity of codimension one Anosov actions of $\RR^k$ on closed manifolds}

\author{Thierry Barbot}
\author{Carlos Maquera}
\thanks{The second author would like to thank CNPq of Brazil for financial
support Grant 200464/2007-8}

\keywords{Anosov action, compact orbit, closing lemma, robust
transitivity, Anosov flow.}

\subjclass[2000]{Primary: 37C85}

\date{\today}

\address{
Thierry Barbot\\
CNRS, UMR 5669, UMPA, ENS Lyon 46, allée d'Italie 69364 Lyon\\
}
\address{
Carlos Maquera\\
CNRS, UMR 5669, UMPA, ENS Lyon 46, allée d'Italie 69364
Lyon\\
and \\
Universidade de S{\~a}o Paulo - S{\~a}o Carlos \\Instituto de
ci{\^e}ncias
matem{\'a}ticas e de Computa\c{c}{\~a}o\\
Av. do Trabalhador S{\~a}o-Carlense 400 \\
13560-970 S{\~a}o Carlos, SP\\
Brazil}

 \email{barbot@umpa.ens-lyon.fr}
 \email{camaquer@umpa.ens-lyon.fr}

 \begin{abstract}
 In this paper, we consider Anosov actions of
 $\RR^k,\ k\geq 2,$ on a closed connected orientable  manifold
 $M$, of codimension one, i.e. such that 
 the unstable foliation associated to some element of $\RR^k$ has dimension one.
 We prove that if the ambient manifold has dimension greater
 than $k+2$, then the action is topologically transitive.
 This generalizes a result of Verjovsky for codimension one Anosov
 flows.
 \end{abstract}

 \maketitle

 \medskip
 \medskip

 \thispagestyle{empty}

%%%%%%%%%%%%%%%%%%%%%%%%%%%%%%%%%%%%%%%%
%%%%%%%%%%%%%%%%%%%%%%%%%%%%%%%%%%%%%%%%
%\tableofcontents

 \section{Introduction}
 It is nowadays a common sense that the Anosov systems
 lie in the central heart of the theory of dynamical systems,
 as the most perfect kind of global hyperbolic behavior.
 It has strong connections with algebra, natural
 examples arising from number field theory or Lie groups theory (see also
 \cite{labourieAnosov} for an example illustrating the deep interplay between Anosov systems and
 representation theory), and also with topology, the dynamics of an Anosov system usually reflecting
 the ambient manifold topology.

 The notion has been introduced by V.V. Anosov in the 60's in \cite{anosov},
 but one should also consider previous works by precursors, including Hadamard,
 Morse, etc...

 An Anosov system is \textit{(topologically) transitive} if it admits a dense orbit.
 There is a quite extensive literature devoted to transitivity for certain classes
 of Anosov systems.
 In particular, by a celebrated result of Newhouse \cite{newhouse} and Franks \cite{franks},
 every codimension one Anosov diffeomorphism on a compact manifold
 is topologically mixing (more than transitive). As a corollary
 from this theorem and \cite[Corollary (6.4)]{franks}, up to finite coverings,
 Anosov diffeomorphisms of codimension one (i.e. such that the unstable 
 subbundle has dimension one) on closed manifolds of dimension $\geq 3$ are topologically
 conjugate to hyperbolic toral automorphisms. For flows, in the
 three-dimensional case, Franks--Williams \cite{FrWi} construct an Anosov flow
 that is not topologically transitive. In the higher dimensional case, Verjovsky
 \cite{verjovsky} proved that codimension one Anosov flows
 on manifolds of dimension greater than three are transitive.

 A natural question arises:
``what about transitivity for actions of higher dimensional groups
 (particularly $\RR^k,\ k\geq 2$)?"

 The development which concerns us here deals with the case of
 Anosov actions of the group $\RR^k$ (some element $r\in \RR^k$ acts normally
  hyperbolically with respect to the orbit foliation).
 This concept was originally introduced by Pugh--Shub \cite{push} in the early
 seventies, and more recently received a strong impetus under the contribution
 of A. Katok and R.J. Spatzier. The rigidity aspects of these actions
 receives nowadays a lot of attention, in the framework of Zimmer program.

 In this paper we undertake the study of transitivity of codimension one Anosov
 actions of $\RR^k, k>1$. An action $\mathbb{R}^{k}$ is called \textit{topologically transitive}
 if it admits a dense orbit. Our main result is the following
 theorem.

 \begin{theorem}
 \label{thm:main}
 Every codimension one Anosov action of $\RR^k$ on a closed manifold
 of dimension greater than $k+2$ is topologically transitive.
 \end{theorem}
 
 The denomination, coming from the usual one for Anosov diffeomorphisms or flows,
 may be confusing: here, codimension one does not mean that 
 the orbits of $\RR^k$ have codimension one, but that the unstable foliation of some
 element of $\RR^k$ has dimension one. See \S~\ref{sub.def}. 

 Note that if a closed $n$-manifold support a codimension one Anosov action of
 $\RR^k$ and $m< k+3,$ then $m=k+2.$ In this case, the Theorem does not hold:
 take the product (cf. Example~\ref{ex:produitflot})
 of the by Franks--Williams example (\cite{FrWi}) by a flat torus
 is a non transitive codimension one Anosov action of
 $\RR^k$ on a $(k+2)$-manifold of the form $N^3\times
 \mathbb{T}^{k-1}$, where $N^3$ is closed three manifold.

 Actually, we will prove slightly more. The theorem above states that under the hypothesis
 there is a dense $\RR^k$-orbit, but we can wonder if there is a one parameter subgroup
 of $\RR^k$ whose orbit on $M$ is dense. Actually, this stronger statement
 does not hold in general: just consider once more as above the product of an Anosov flow, transitive or not,
 by a flat torus. However, it is nearly true, in a weak sense,
 as explained just below.

 An element of $\RR^k$ is said Anosov if it acts normally
 hyperbolically with respect to the orbit foliation. Every connected component of
 the set of Anosov elements is an open convex cone in $\RR^k$, called a chamber. More generally,
 a \textit{regular subcone} $\cC$ is an open convex cone in $\RR^k$ containing only Anosov elements.
 One should consider $\cC$ as a semi-group in $\RR^k$: the sum of two elements in the
 cone still lies in the cone. The $\cC$-orbit of a point $x$ in $M$ is the subset
 comprising the iterates $\phi^a(x)$ for $a$ describing $\cC$.

 \begin{theorem}
 \label{thm:main2}
 Let $\phi$ be a codimension one Anosov action of $\RR^k$ on a closed manifold $M$
 of dimension greater than $k+2$. Then any regular subcone $\cC$ admits a dense orbit in
 $M$.
 \end{theorem}

 Theorem~\ref{thm:main} is obviously a direct corollary of Theorem~\ref{thm:main2}.
 On the other hand, given $a$ in $\RR^k$, we can apply Theorem~\ref{thm:main2}
 to every small regular cone containing $a$: hence we can loosely have in
 mind that, up to arbitrarily small errors, $\phi^a$ admits a dense orbit.

 The cornerstone of the proof is the study of the codimension foliation $\cF^s$
 tangent to the stable subbundle of Anosov elements in $\cC$. The unstable
 foliation $\cF^{uu}$ for $\cC$ has dimension one, and the first step is to prove that
 every leaf of $\cF^{uu}$ admits an affine structure, preserved by the action of $\RR^k$
 (cf. Theorem~\ref{thm:AfinStruct}). The existence of this affine structure provides
 a very good information about the transverse holonomy of $\cF^s$,
 giving \textit{in fine}, through the classical theory of codimension one foliations,
 many information about the topology of the various foliations involved. In particular,
 the orbit space of the lifting of $\phi$ to the universal covering $\widetilde{M}$
 is a Hausdorff manifold, homeomorphic to $\RR^{n-k}$ (cf. Theorem~\ref{thm:orbitspace}).

 On the other hand, one can produce a generalization for $\RR^k$-actions of the
 classical Spectral decomposition Theorem (Theorem~\ref{thm:espdecomp}). It allows to reduce
 the proof of the transitivity to the proof that stable leaves are dense (Lemma~\ref{lem:transitivity}).
 Now, if some leaf of $\cF^{s}$ is not dense, then there must be some \textit{non bi-homoclinic orbit}
 of $\RR^k$ (cf. Proposition~\ref{prop:pointbihomoclinic}). One then get the final result
 by using some clever argument, involving Jordan-Sch\"{o}nflies Theorem, and already
 used in Verjovsky proof as rewritten in  \cite{barthese} or \cite{matsumoto}.

 Actually, all the strategy above mostly follows the guideline used in Verjovsky proof,
 but is more than a simple transposition. New phenomena arise, even enlightening the case
 of Anosov flows.

 -- \textit{Convex cones:} The case $k=1$ is somewhat greatly simplified
 by the fact that regular subcones in $\RR^k$ are simply half-lines, and that the only
 non Anosov element of $\RR^k$ is the origin $0$. One can compare the classical Closing Lemma
 with the general version, more technical in its statement: Theorem~\ref{thm:ClosingLemma}.

-- \textit{Reducibility:} Given an Anosov $\RR^k$ action on some manifold $M$, one can always
 take the product $M \times \TT^l$ by some torus $\TT^l = \RR^l/\ZZ^l$ and consider the locally
 free action of $\RR^{k+l}$ on this product manifold. Then this action is still Anosov. This construction
 can be generalized to twisted products through a representation $\rho: \pi_1(M) \to \RR^l$
 (cf \S~\ref{ex:bundle}). Of course, this construction gives examples with $k > 1$, hence doesn't
 appear in the case of Anosov flows.
 Therefore, an important step is to put aside these examples. In Theorem~\ref{thm:ReducingAction},
 we prove that every codimension one action of $\RR^k$ splits uniquely as a principal
 torus bundle over some manifold $\bar{M}$ such that the $\RR^k$ actions permutes the fibers,
 and thus induces an action on $\bar{M}$. Moreover, the fibers are precisely the orbits
 of some subgroup $H_0 \subset \RR^k$, and the induced action is Anosov. Finally, this
 splitting is maximal, i.e. $\bar{M}$ cannot be decomposed further: it is \textit{irreducible}.
 Many properties, among them transitivity, is obviously satisfied by the Anosov action on $M$ if and
 only if it is satisfied by the induced action on $\bar{M}$. Therefore, the
 proof of Theorem~\ref{thm:main2} reduces to the irreducible case.

 The irreducibility of an Anosov $\RR^k$ action can be equivalently defined as
 requiring that the codimension one stable foliation $\cF^s$ has trivial holonomy
 cover; in a less pedantic way, it means that homotopically non-trivial loops in leaves have non-trivial
 holonomy (cf. Remark~\ref{rk:defreduc}). Irreducible Anosov actions enjoy many nice topological
 properties. Among them (cf. Proposition~\ref{pro:OrbIrrecAct}):

 \textit{Let $\phi$ be an irreducible codimension one action of $\RR^k$ on a manifold $M$.
 Then the isotropy subgroup of every element of $M$ is either trivial, or a lattice
 in $\RR^k$.}

 We can observe, as a corollary, that if a codimension one Anosov action of $\RR^k$
 admits an orbit homeomorphic to $\TT^{k-1} \times \RR$, then it is a twisted
 product by flat tori $\TT^{k-1}$ over an Anosov flow.

 The paper is organized as follows: in the preliminary section~\ref{sec.preli}, we give definitions, and present
 first results, as the generalized Closing Lemma for actions of $\RR^k$, and the spectral decomposition
 of the non-wandering set as a finite union of basic blocks. In \S~\ref{sec:examples} we present the
 known examples of Anosov actions of codimension one. In \S~\ref{sec:reducing} we establish the reduction
 Theorem~\ref{thm:ReducingAction} (it includes the proof of the $\RR^k$-invariant affine structures
 along unstable leaves). In \S~\ref{sec:proof} we prove the Main Theorem~\ref{thm:main2}.
 In the last section~\ref{sec:conclusion} we give additional comments, and present
 forecoming works in progress.

 \vspace{.5cm}
 \textbf{Acknowledgments.} This paper was written while the second author stayed at Unit\'e
de Math\'ematiques Pures et Appliqu\'ees, \'Ecole Normale Sup\'erieure de
Lyon. He thanks the members of UMPA, especially Professor Etienne
Ghys for his hospitality.

\section{Preliminaries}
\label{sec.preli}

\subsection{Definitions and notations}
 Now, we outline some basic results about actions of $\RR^k$ which
 will be used in the proof of the  main theorem. Recall that, for any
 action $\phi :\RR^k \times M \to M$ of $\RR^k$ on a manifold $M,$
 $\mathcal{O}_p:=\{\phi(\om, p),\om \in \RR^2 \}$ is the orbit of
 $p \in M$ and $\Gamma_p := \{ \om \in\RR^k : \phi(\om, p) = p \}$ is
 called the isotropy group of $p.$ The action $\phi$
 is said to be \textit{locally free} if the isotropy group of every point is
 discrete. In this case the orbits are diffeomorphic to
 $\mathbb{R}^{\ell} \times \mathbb{T}^{k-\ell}$, where $0\leq \ell \leq k$.

Let $\mathcal{F}$ be a continuous foliation on a manifold $M$. We
denote the leaf that contains $p \in M$ by $\cF(p)$. For an open
subset $U$ of $M$,
 let $\cF|_U$ be the foliation on $U$ such that
 $(\cF|_U)(p)$ is the connected component of $\cF(p) \cap U$
 containing $p \in M$.
A coordinate $\vphi=(x_1,\cdots, x_n)$ on $U$
 is called {\it a foliation coordinate} of $\cF$
 if $x_{m+1},\cdots,x_n$ are constant functions
 on each leaf of $\cF|_U$, where $m$ is the dimension of $\cF$.
A foliation is of class $C^{r+}$ is if it is covered by
 $C^{r+}$ foliation coordinates.
We denote the tangent bundle of $M$ by $TM$. If $\cF$ is a $C^1$
foliation, then we denote the tangent bundle
 of $\cF$ by $T\cF$.
%For a homeomorphism $h$ from $M$ to another manifold $M'$,
% we define a foliation $H(\cF)$ on $M'$
% by $H(\cF)(p)=H(\cF(H^{-1}(p)))$.

 We fix a Riemannian metric $\|$, and denote by
 $d$ the associated distance map on $M$.

%%%%%%%%%%%%%%%%%%%%%%%%%%%%%%%%%%%%%%%%%%
 \subsection{Anosov $\RR^k$-actions}
 \label{sub.def}
%%%%%%%%%%%%%%%%%%%%%%%%%%%%%%%%%%%%%%%%%%%
 Let us recall the definitions and basic properties of Anosov
 actions.

 \begin{definition}
 {\rm
 Let $M$ be a $C^\infty$ manifold and $\phi$ a locally free $C^{1+}$ action of
 $\RR^k$ on $M$. By $T\phi$, we denote the $k$-dimensional subbundle of $TM$
 that is tangent to the orbits of $\phi$.
 \begin{enumerate}
   \item We say that $a\in \RR^k$  is an \textit{Anosov element} for $\phi$ if $g=\phi(a,\cdot)$
   acts normally hyperbolically with respect to the orbit foliation. That is, there exist real
  numbers $\lambda > 0,\ C > 0$ and a continuous $Dg$-invariant splitting of the tangent bundle
  $$
  TM=E_a^{ss}\oplus T\phi\oplus E_a^{uu}
  $$
  such that
  $$
  \begin{array}{cc}
   \|Dg^n|_{E_a^{ss}}\|\leq Ce^{-\lambda n} & \forall n>0 \\
   \|Dg^n|_{E_a^{uu}}\|\leq Ce^{\lambda n} & \forall n<0
  \end{array}
  $$
   \item Call $\phi$ an \textit{Anosov action} if some $a\in \mathbb{R}^k$ is
   an Anosov element for $\phi$.
  \end{enumerate}
 }
 \end{definition}

 Hirsch, Pugh and Shub developed the basic theory of normally hyperbolic transformations
 in \cite{hirpush}. As consequence of this we obtain that the splitting is H\"older continuous
 and the subbundles $E_a^{ss}$, $E_a^{uu}$, $T\phi \oplus E_a^{ss}$, $T\phi \oplus E_a^{uu}$
 are integrable.
 The corresponding foliations, $\cF_a^{ss}$, $\cF_a^{uu}$, $\cF_a^{s}$, $\cF_a^{u}$, are called
 \textit{the strong stable foliation, the strong unstable foliation, the weak stable foliation},
 and \textit{the weak unstable foliation}, respectively.

 From now $\phi$ is an Anosov action of $\mathbb{R}^k$ on $M$,
 and $a$ an Anosov element fixed once for all.
 For simplicity, the foliations corresponding to $a$ will be denoted by
 $\cF^{ss}$, $\cF^{uu}$, $\cF^{s}$ and  $\cF^{u}$. For all $\delta >0$, $\cF_{\delta}^{i}(x)$
 denote the open ball in $\cF^i(x)$ under the induced metric which centering at $x$ with radius
 $\delta$, where $i=ss,uu,s,u.$

 \begin{theorem}[of product neighborhoods]
 \label{thm:local product}
 Let $\phi:\RR^k \times M \to M$ be an Anosov action. There exists a
 $\delta_0 >0$ such that for all $\delta \in (0,\delta_0)$ and for
 all $x\in M,$ the applications
 $$
 [\cdot,\cdot]^u:\cF^s(x)\times \cF^{uu}(x)\to M; \ \
 [y,z]^u=\cF^s_{2\delta}(z)\cap \cF^{uu}_{2\delta}(y)
 $$
 $$
 [\cdot,\cdot]^s:\cF^{ss}(x)\times \cF^{u}(x)\to M; \ \
 [y,z]^s=\cF^{ss}_{2\delta}(z)\cap \cF^{u}_{2\delta}(y)
 $$

 are homeomorphisms on their images.
 \end{theorem}

 \begin{remark}
 \label{rk.zero}
 {\rm
Every foliation $\cF^{ss}_a$, $\cF^{uu}_a$, $\cF^{s}_a$ or  $\cF^{u}_a$ is preserved by every diffeomorphism commuting with $a$. In particular, it is
 $\RR^k$-invariant.  Another standard observation is that, since every compact domain in a leaf of $\cF_a^{ss}$ (respectively
 of $\cF_a^{uu}$) shrinks to a point under positive (respectively negative) iteration by $\phi^a$, every leaf
 of $\cF_a^{ss}$ or $\cF_a^{uu}$ is a \textit{plane}, i.e. diffeomorphic to $\RR^\ell$ for some $\ell$.

 Let $F$ be a weak leaf, let say a weak stable leaf. For every strong stable leaf $L$ in $F$, let
 $\Gamma_L$ be the subgroup of $\RR^k$ comprising elements $a$ such that $\phi^a(L)=L$,
 and let $\cO_L$ be the saturation of $L$ under $\phi$. Thanks to Theorem~\ref{thm:local product} we have:
 \begin{itemize}
 \item $\cO_L$ is open in $F$,
 \item $\Gamma_L$ is discrete.
 \end{itemize}
 Since $F$ is connected, the first item implies $F=\cO_L$: the $\phi$-saturation of a strong leaf  is an entire weak leaf.
 Therefore, $\Gamma_L$ does not depend on $L$, only on $F$.
 The second item implies that the quotient $P=\Gamma_L\backslash\RR^k$ is a manifold, more precisely,
 a flat cylinder, diffeomorphic to $\RR^p\times\TT^q$ for some $p$, $q$.
 For every $x$ in $F$, define $p_F(x)$ as the equivalence class $a+\Gamma_L$ such that
 $x$ belongs to $\phi^a(L)$.  The map $p_F: F \to P$ is a locally trivial fibration
 and the restriction of $p_F$ to any $\phi$-orbit in $F$ is a covering map.
 Since the fibers are contractible (they are leaves of $\cF^{ss}$, hence planes), the fundamental group
 of $F$ is the fundamental group of $P$, i.e. $\Gamma_L$ for any strong stable leaf $L$ inside $F$.

 Observe that if $\cF^{ss}$ is oriented, then the fibration $p_F$ is trivial: in particular, $F$ is diffeomorphic to
 $P \times \RR^{p}$, where $p$ is the dimension of $\cF^{ss}$.

 Of course, analogous statements hold for the strong and weak unstable leaves.
 }
 \end{remark}

 \vspace{.5cm}
  We say $\phi$ is a {\it codimension-one} Anosov action
 if $E^{uu}_a$ is one-dimensional for some $a$ in $\mathbb{R}^k$. In this case,
 we will always assume that the fixed Anosov element has one dimensional strong stable foliation.

 \begin{remark}
 \label{rk.first}
 {\rm
 Let $\mathcal{A}=\mathcal{A(\phi)}$ be the set of Anosov elements of $\phi$.
 }
 \begin{enumerate}
   \item $\mathcal{A}$ is always an open subset of $\mathbb{R}^k$.
   {\rm
   In fact, by the structural stability theorem
   for normally hyperbolic transformations by Hirsch, Pugh and Shub a map $C^1$-close to a normally
   hyperbolic transformation is again normally hyperbolic for a suitable foliation \cite{hirpush}. For an
  element in $\mathbb{R}^k$ close to an Anosov element, this suitable foliation is
  forced to be the orbit foliation of the action.
  }
   \item Every connected component of $\mathcal{A}$ is an open convex cone in
   $\mathbb{R}^k$.
   {\rm
   Let $a$ be an Anosov element. Every element near $a$ must share the same stable and unstable
   bundles, therefore, all the Anosov elements in the same connected component than $a$
   admits the same stable/unstable splitting. The contracting or expanding property along
   a given bundle is stable by composition and by multiplication of the generating vector field
   by a positive constant factor; it follows that the connected component is
   a convex cone, as claimed.

   We call such a connected component a \textit{chamber}, by analogy with the case of
   Cartan actions. More generally, a \textit{regular subcone} is an open convex cone
   contained in a chamber.

   }
   \item If $\mathcal{A}_a$ is the chamber containing $a,$ then
   $\cF_a^{ss}=\cF_b^{ss}$, $\cF_a^{uu}=\cF_b^{uu}$, $\cF_a^{s}=\cF_b^{s}$
   and
   $\cF_a^{u}=\cF_b^{u}$ for all $b\in \mathcal{A}_a.$
   \item Any $\phi$-orbit whose isotropy subgroup
   contains an Anosov element $v$ is compact.
   {\rm
   Indeed, let $y$ be a point in the closure of the orbit. It is clearly fixed by $\phi^v$, and there is a  local cross-section
   $\Sigma$  to $\phi$ containing $x$ such that for every $z$ in $\Sigma$ near $y$ the image $\phi^v(z)$ lies in $\Sigma$.
   Then, $y$ is a fixed point of $\phi^v$ of saddle type, in particular, it is an isolated $\phi^v$-fixed point. Hence
   $y$ lies in the $\phi$-orbit of $x$.
    }
   \end{enumerate}
 \end{remark}

  Another standard fact about Anosov $\mathbb{R}^k$-actions is an
 Anosov-type closing lemma which is a straightforward generalization
 of a similar statement for Anosov flows (cf. \cite[Theorem 2.4]{KatSpatz1}).

 \begin{theorem}[Closing Lemma]
 \label{thm:ClosingLemma}
 Let $a\in \mathbb{R}^k$ be an Anosov element of an Anosov $\mathbb{R}^k$-action
 $\phi$ on a closed manifold $M$. There exist positive constants $\varepsilon_0$, $C$
 and $\lambda$ depending continuously on $\phi$ in the $C^1$-topology and $a$ such
 that: if for some $x\in M$ and $t\in \mathbb{R}$
 $$d(\phi(ta,x), x) < \varepsilon_0,$$
  then there exists a point
  $y\in M$, a differentiable map $\gamma:[0, t]\to \mathbb{R}^k$
  such that for all $s\in [0, t]$ we have
  \begin{enumerate}
    \item $d(\phi(sa,x),\phi(\gamma(s),y))<Ce^{-\lambda(\min\{s,t-s\})}
         d(\phi(ta,x), x);$
    \item $\phi(\gamma(t),y)=\phi(\delta,y)$ where $\|\delta\| < C d(\phi(ta,x),
    x);$
    \item $\|\gamma'-a\|< Cd(\phi(ta,x), x)$.
  \end{enumerate}
   \end{theorem}

 \begin{remark}
 \label{rk.second}
 {\rm
 Let $\cC$ be a regular subcone containing $a$ (for example, a chamber). Once $a$ is fixed, item (3) in the Theorem above
 implies that if $d(\phi(ta,x), x)$ is sufficiently small, the velocity $\gamma'$ lies in
 $\cC$, therefore, that the image of $\gamma$ is contained in $\cC$.  Moreover, once
 more if $d(\phi(ta,x), x)$ is sufficiently small, item (2) implies that $\gamma(t)-\delta$ belongs to
 $\cC$. According to Remark~\ref{rk.first} the orbit of $y$ is compact.
 }
  \end{remark}

 \begin{definition}[The nonwandering set]
 \label{def:nonwandering}
 {\rm
 A point $x\in M$ is \textit{nonwandering} with respect to a regular subcone $\cC$
 if for any open set $U$ containing $x$ there is a $v\in \cC, \ \|v\|>1,$
 such that $\phi^v(U)\cap U \neq \emptyset,$ where
 $\phi^v=\phi(v,\cdot)$.
 The set of all nonwandering points, with respect to $\cC$, is denoted by
 $\Omega(\cC)$.
 }
 \end{definition}

 By using the Closing Lemma for Anosov $\mathbb{R}^k$-actions we
 obtain:

  \begin{proposition}
 \label{prop:closCompactO}
  For any regular subcone $\cC$, the union of compact orbits of $\RR^n$ is dense in $\Omega(\cC)$.
 \end{proposition}
 \begin{proof}
 For $x\in \Omega(\cC)$ and $\varepsilon > 0$ denote by
 $U_{\varepsilon}$ the $\varepsilon /(2C+1)$-neighborhood of $x$ in
 $M$, where $C$ is as in the Closing Lemma. Then there exists $v\in \cC$
  such that $\phi^v(U_{\varepsilon})\cap U_{\varepsilon}\neq \emptyset.$
 For $y\in \phi^{-v}(U_{\varepsilon})\cap U_{\varepsilon}\neq \emptyset$ we have
 $d(\phi^v(y),y)< 2\varepsilon /(2C+1)$ and hence by the Closing
 Lemma and Remark~\ref{rk.second} there is a point $z$ such that $\phi^v(z)=z, \ \mathcal{O}_z$ is compact and
 $d(y,z)<Cd(\phi^v(y),y)$, consequently
 $d(z,x)\leq d(y,x)+d(y,z)<\varepsilon.$ It prove that ${\rm Comp}(\phi)$
is dense in $\Omega(\cC)$ and finishes the proof.
 \end{proof}

 \begin{remark}
 \label{rk.nonwanderingflow}
 Let $a$ be any non-trivial element of $\RR^k$. The nonwandering set $\Omega(\phi^{ta})$
 of the (semi-)flow generated by $a$ is clearly contained in $\Omega(\phi)$. On the other
 hand, the nonwandering set of any linear flow on a torus is the entire torus. Hence,
 compact orbits of $\RR^k$ are contained in $\Omega(\phi^{ta})$. Hence,
 it follows from the proposition above that the nonwandering sets $\Omega(\phi^{ta})$ and
 $\Omega(\phi)$ coincide. In particular, the nonwandering set $\Omega(\mathcal{C})$ is independent
 from the regular subcone $\cC$.
 \end{remark}

 \begin{lemma}
 \label{le:compactcausal}
 The isotropy subgroup of any compact orbit contains an element in $\cC$.
 \end{lemma}
\begin{proof}
Let $R>0$ such that every Euclidean ball of radius $R$ in $\RR^k$ intersects every orbit of the isotropy
subgroup $\Gamma$. Let $B$ be a closed Euclidean ball of radius $r$  in the open convex cone $\cC$.
Then, for $t>R/r$, the ball $tB$ is contained in $\cC$ and has radius $>R$, and thus intersects
the $\Gamma$-orbit of $0$. The lemma follows.
\end{proof}

The Riemannian metric induces an area form on every $\phi$-orbit.

\begin{lemma}
\label{le:raghu}
For every $C>0$, there is only a finite number of compact $\phi$-orbits of area $\leq C$.
\end{lemma}

\begin{proof}
Assume by contradiction the existence of an infinite sequence of
distinct compact orbits $\cO_n$ of area $\leq C$. For each of them,
let $\Gamma_n$ be the isotropy group of $\cO_n$: it is an element of
$\cR=\operatorname{GL}(k,\RR)/\operatorname{SL}(k,\ZZ)$, the space
of lattices in $\RR^k$. Since $\phi$ is locally free, the length of
elements of $\Gamma_n$ is uniformly bounded from below,
independently from $n$. By the Mahler's criterion (\cite{raghu}), it
ensures that, up to a subsequence, the $\Gamma_n$ converges to some
lattice $\Gamma_\infty$. In particular, for every $v_\infty$ in
$\Gamma_\infty$, there is a sequence of elements $v_n$ of $\Gamma_n$
converging in $\RR^k$ to $v_\infty$. Furthermore, according to (the
proof of) Lemma~\ref{le:compactcausal}, we can select $v_\infty$ in
$\cC$. Up to a subsequence, we can also pick up a sequence of
elements $x_n$ in each $\cO_n$ converging to some $x_\infty$ in $M$.
Then, since $\phi^{v_n}(x_n)=x_n$, at the limit we have
$\phi^{v_\infty}(x_\infty)=x_\infty$. Since $v_\infty$ is Anosov,
the $\phi$-orbit $\cO_\infty$ of $x_\infty$ is compact. Consider a
local section $\Sigma$ to $\phi$ containing $x_\infty$: the first
return map on $\Sigma$ along the orbit of $\phi^{v_\infty}$ is
hyperbolic, admitting $x_\infty$ as an isolated fixed point. On the
other hand, by pushing slightly along $\phi$, we can assume without
loss of generality that every $x_n$ belongs to $\Sigma$. Since the
$v_n$ converges to $v_\infty$, the $\phi^{v_n}$-orbit of $x_n$
approximates the $\phi^{v_\infty}$-orbit of $x_\infty$, showing that
the $x_n$ are also fixed points of the first return map. It is a
contradiction, since they accumulate to the isolated fixed point
$x_\infty$.
\end{proof}

 \begin{theorem}[Spectral decomposition]
 \label{thm:espdecomp}
 Let $M$ be a closed smooth manifold and let $\phi$ be an Anosov action on $M$.
 The nonwandering set of $\phi$ can be partitioned into a finite number of
 $\phi$-invariant closed subsets, called basic blocks:
 $$
 \Omega = \bigcup_{i=1}^{\ell} \Lambda_i
 $$
 such that for every regular subcone $\cC$, every $\Lambda_i$ is $\cC$-transitive, i.e. contains
 a dense $\cC$-orbit.
 \end{theorem}
 \begin{proof}
 Let ${\rm Comp}(\phi)$ be the set of compact orbits of $\phi$. By Proposition \ref{prop:closCompactO}
 we have $\overline{{\rm Comp}(\phi)}=\Omega(\phi).$
  We define a relation on ${\rm Comp}(\phi)$ by: $x\sim y$ if and
  only if $\cF^u(x)\cap \cF^s(y)\neq \emptyset$ and $\cF^s(x)\cap \cF^u(y)\neq \emptyset$
  with both intersections transverse in at least one point. We want
  to show that this is an equivalence relation and obtain each
  $\Lambda_i$ as the closure of an equivalence class.

  Note that $\sim $ is trivially reflexive and symmetric. In order to check
  the transitivity suppose that $x,y,z\in {\rm Comp}(\phi)$ and
  $p\in \cF^u(x)\cap \cF^s(y), $ $q\in \cF^u(y)\cap \cF^s(z)$ are
  transverse intersection points. There exists $v\in \mathcal{A}_a$ such that
  $\phi^v(x)=x$. Since the images of a ball around $p$ in
  $\cF^u(p)=\cF^u(x)=\phi^v(\cF^u(x))$ accumulate on $\cF^u(y)$, we
 obtain that $\cF^u(x)$ and $\cF^s(z)$ have a transverse
 intersection. Analogously, we obtain that $\cF^s(x)$ and $\cF^u(z)$ have a transverse
  intersection.

 By Theorem \ref{thm:local product} any two sufficiently near points
 are equivalent, so by compactness we have finitely many equivalence
 classes whose (pairwise disjoint) closures we denote by
 $\Lambda_1,\Lambda_2,\dots, \Lambda_{\ell}$.

 It remains to show that every $\Lambda_i$ is $\cC$-transitive for every regular subcone $\cC.$
 Notice first that if $p\in \Lambda_i \cap {\rm Comp}(\phi)$ and $p\sim q$
 with $q\in {\rm Comp}(\phi)$, then there is $z\in \cF^u(p)\cap \cF^s(q).$
 Let $v\in \cA_a$ such that $\phi^v(p)=p$.
 As the iterates under $\phi^v$
 of a ball around $z$ in $\cF^u(p)=\cF^u(z)=\phi^v(\cF^u(z))$ accumulate on $\cF^u(q)$,
 and since $z$ belongs to $\cF^s(q)$, we
 obtain that $\cF^u(p)$ is dense in $\Lambda_i \cap
 {\rm Comp}(\phi)$, hence in $\Lambda_i.$

 Now, for the transitivity, we need to show that for any two open sets
 $U$ and $V$ in $\Lambda_i$ there exists $v\in \cC$ such
 that $\phi^v(U)\cap V\neq \emptyset.$ The density of compact orbits
 in $\Lambda_i$ implies the existence of $p\in U$ and $v\in \cC$
 such that $\phi^v(p)=p$  (cf. Lemma~\ref{le:compactcausal}).
 Let $\cF^{uu}_{\delta}(p)$ be a
 neighborhood of $p$ in $\cF^{uu}(p)$ that is contained in $U.$
 Since the $\RR^k$-orbit of $p$ is compact, there is a compact domain $K$
 in $\RR^k$ so that the leaf
 $\cF^{u}(p)$ is equal to $K\cdot \cup_{j=0}^{\infty}\phi^{jv}(\cF^{uu}_{\delta}(p))$.
 Since this leaf is dense in $\Lambda_i$, there exists for every $m \in \NN$ sufficiently big
 an element $g_m\in K$
 such that $V\cap [g_m \cdot \cup_{j=0}^{m}\phi^{jv}(\cF^{uu}_{\delta}(p))]\neq \emptyset,$
 hence, $V\cap [\cup_{j=0}^{m}\phi^{jv+g_m}(\cF^{uu}_{\delta}(p))]\neq \emptyset.$
 Thus $V\cap \phi^{mv+g_m}(\cF^{uu}_{\delta}(p))\neq \emptyset,$
 consequently, $V\cap \phi^{mv+g_m}(U)\neq \emptyset.$ The Theorem follows, since
 for $m$ sufficiently large, $mv+g_m = m(v+g_m/m)$ lies in $\cC$.
 \end{proof}

\section{Examples}
\label{sec:examples}

Let us give some examples of Anosov actions of
$\RR^k$. We will especially focus on codimension one examples.

 \begin{example}
{\rm Let $G$ be a real semi-simple Lie group, with Lie algebra
$\cG$, $\Gamma$ a torsion-free uniform lattice in $G$, and $A$ a
split Cartan subgroup of $G$. The centralizer of $A$ in $G$ is a
product $AK$, where $A$ commutes with $K$. Then the action at the
right induces a $\RR^k$-action on the compact quotient $M =
\Gamma\backslash{G}/K$. An essential starting point in the theory of
root systems has a strong dynamical system flavor: this action is
Anosov! More precisely, the classical first step is to prove that
the adjoint action of $A$ on $\cG$ preserves a splitting~:

$$
\cG = \mathcal{K} + \mathcal{A} + \sum_{\alpha \in \Sigma} \mathcal{G}^\alpha
$$

where $\mathcal{K}$, $\mathcal{A}$ are the Lie algebras of $K$, $A$,
and where every $\alpha$ (the {\em roots}) are linear forms
describing the restriction of the adjoint action of $a$ on
$\cG^\alpha$: it is simply the multiplication by $\alpha(a)$. The
classical way is then to prove that the elements $a$ of $A$ for
which $\alpha(a) \neq 0$ is a Zariski open subset, and these
elements, called {\em regular}, are precisely the ones which are
Anosov in our terminology for $\RR^k$-action. They form an union of
open convex cones, called Weyl chambers,  of the form $\{ \alpha >
0; \forall \alpha \in \Sigma^\ast \}$ where $\Sigma^\ast$ is a
subsystem of a certain kind, called {\em reduced root system.}

This family of examples, called \textit{Weyl chamber flows} in \cite{KatSpatz1},
is certainly one the the most interesting, but is never of codimension one,
except in the case $G = \operatorname{SL}(2,\RR)$. Indeed, the root system is always
equal to its own opposite. Hence if the associated $\RR^k$ Anosov action has
codimension one, then $\Sigma$ contains exactly two elements, and our assertion follows.
In this very special case, the examples we obtain are
Anosov flows (i.e. $k=1$), and more precisely, up to finite coverings, geodesic
flows of compact Riemannian
surfaces with constant curvature $-1$.

This dynamical feature in algebra is very useful. For example, the (simple) fact that
Anosov actions admit compact orbits implies that every uniform lattice $\Gamma$
in $G$ admits a conjugate $g\Gamma g^{-1}$ which is a lattice in $A$. In particular,
$\Gamma$ contains a free abelian subgroup of the same rank than $G$.
}
\end{example}

\begin{example}
\label{ex:suspension} {\rm Consider an action of $\ZZ^k$ on a closed
manifold $S$. The \textit{suspension} of this action is the quotient
$M$ of $S \times \RR^k$ by the relation identifying each $(x, u)$
with $(k.x, u+k)$ for every $k$ in $\ZZ^k$. The translation on the
second factor $\RR^k$ induces a $\RR^k$-action on $M$. It is easy to
prove that this action is Anosov if and only some element of $\ZZ^k$
induces an Anosov diffeomorphism on $S$. Observe also that this
action has codimension one if and only if one of the Anosov element
of $\ZZ^k$ has codimension one. Hence, by Franks-Newhouse Theorem
reported in the introduction, if the suspension has codimension one,
then codimension one Anosov elements of $\ZZ^k$ are (up to finite
coverings and topological conjugation) hyperbolic toral automorphism
on some torus $\TT^n$. Every homeomorphism of the torus commuting
with a hyperbolic toral automorphism is also an automorphism (i.e.
linear). Hence the only possible examples of codimension one
suspensions are the ones described below, arising from number field
theory, maybe after restriction to a subgroup of $\ZZ^k \subset
\operatorname{Aut}(\TT^n)$.

The suspension process can be generalized to a version including Weyl chamber flows
(see twisted Weyl chamber flows in \cite{KatSpatz1, KatSpatz2}), but this new family of examples
are never of codimension one.

}
\end{example}

\begin{example}
\label{ex:diffeo}
{\rm
Generically, the centralizer of an Anosov diffeomorphism $f$ reduces to the iterates $f^k \;\; (k \in \ZZ)$
(see \cite{YoccozPalis}). Hence the construction of Anosov actions of $\ZZ^k$ for $k \geq 2$ requires special features.

Let $K=\QQ[\alpha]$ be a field extension of the field $\QQ$ of finite degree $n$,
$\cO_K$ the ring of algebraic integers of $K$, and $\cO^\ast_K$ the group of units
of $\cO_K$. Then, the quotient of $K \otimes \RR$ by the additive action of
$\cO_K$ is a compact torus of dimension $n$, on which $\cO^\ast_K$ acts by multiplication.
According to Dirichlet unit Theorem, the torsion-free part of $\cO^\ast_K$
is isomorphic to $\ZZ^k$, with $k=r_1+r_2-1$ where $r_1$ is the number of real embeddings
and $r_2$ the number of conjugate pairs of complex embeddings of $K$. Hence every
finite extension of $\QQ$ naturally provides an action of $\ZZ^k$ on a torus.
More precisely, the real and complex embeddings provide altogether a realization
of $K \otimes \RR$ as a vector subspace of $\RR^{r_1}\oplus\CC^{r_2}$, preserved by the
multiplicative action of $\cO^\ast_K$, which is diagonalizable, the eigenvalues being
the various conjugates. Hence, this action is Anosov
if and only if some unit has no conjugate of norm $1$.

A concrete way to produce such examples is to take some algebraic number $\alpha$ admitting
no conjugate of norm $1$, and to consider the extension $K=\QQ[\alpha]$. Of course, one can
forget part of the unit group and just consider some subgroup. It is actually what we do when
defining linear Anosov diffeomorphisms.

In order to get codimension one $\ZZ^k$-actions, it is sufficient to select as algebraic integer $\alpha$
any \textit{Pisot number}, which is, by definition precisely a real algebraic integer $\alpha$
exceeding $1$, and such that its conjugate elements are all less than $1$ in absolute value.
Concretely, examples of Pisot numbers are roots of $x^3-x-1$, $x^4-x^3-1$,
etc...

}
\end{example}

 \begin{example}
 \label{ex:produitflot}
 {\rm Let $N$ be a $n$-dimensional manifold supporting a codimension one Anosov flow
 (clearly, $n\geq 3$).
 We construct a codimension one action of $\RR^k,\ k\geq 2,$ on $M = N \times \TT^{k-1}$. Consider
 the coordinate system $(x, \th)$ in $M, x \in N, \th \in \TT^{k-1}.$ In what follows,
 for a real function
 $a(x, \th)$,  by  $ a(x, \th)\frac{\partial}{\partial x} $ we mean
 $ a_1 \frac{\partial}{\partial x_1} +\cdots + a_n \frac{\partial}{\partial x_n}
 $ where $x_1,\dots, x_n$ are coordinates in $N.$

 Let $\phi \in A^1(\RR^k, M)$ be defined by $X_1$ and $Y_1,\dots Y_{k-1}, $  such that
 $X_1 = a(x) \frac{\partial}{\partial x}$ is a codimension one Anosov flow
 in $N$ and $X_j:= \frac{\partial}{\partial \th_j},$ where $\theta_1,\dots, \theta_{k-1}$
 are coordinates in $\TT^{k-1}.$ Then
 \begin{itemize}
   \item $\phi$ is a codimension one Anosov action of $\RR^k$ on $M$.
   \item for $n>3,$ by Verjovsky Theorem, $\phi$ is transitive.
   \item if $n=3$ and $X_1$ is the Anosov flow defined by
   Franks-Williams in \cite{FrWi}, then $\phi$ is a codimension one
   Anosov action of $\RR^k$ on the $(k+2)$-manifold $M$ which is not
   transitive.
 \end{itemize}
 }
 \end{example}

 \begin{example}
 \label{ex:bundle}
{\rm
 More generally, let $\phi$ be an Anosov $\RR^k$-action on a closed $n$-dimensional manifold $M$,
 and let $p: \widehat{M} \to M$ be a principal flat $\TT^\ell$-bundle over $M$. By flat,
 we mean that it is equipped with a flat $\TT^\ell$-invariant connection, i.e. a $n$-dimensional foliation
 $\cH$ transverse to the fibers of $p$ and preserved by $\TT^\ell$ (there is a 1-1 correspondence between
 principal flat $\TT^k$-bundles and group homomorphisms $\rho: \pi_1(M) \to \TT^k$).
Then, the $\RR^\ell$-action on $M$ lifts uniquely as a $\RR^\ell$ action tangent to $\cH$. Moreover,
this action commutes with the right action of $\TT^\ell$ tangent to the fibers. Hence, both action
combine to a $\RR^{k+l}$-action on $\widehat{M}$, which is clearly Anosov.
}

  \end{example}

 \section{Reducing codimension one Anosov actions}
 \label{sec:reducing}
 In this section we to show that any codimension one Anosov action, up to a reduction through
 a principal torus bundle, has the several topological properties, including:

 \begin{itemize}
 \item The universal covering of the ambient manifold is diffeomorphic ro $\RR^n$;
 \item The fundamental group of every compact orbit injects into the fundamental group of the ambient
 manifold;
 \item The holonomy of the codimension one foliation along a homotopically non-trivial
 loop in a leaf is non-trivial.
 \end{itemize}

 The crucial ingredient is the construction along each strong leaf of dimension one of
 an affine structure, preserved by the action of $\phi$.

 \subsection{Affine structures over the leaves strong unstables}

 We begin by remembering that an affine structure of class $C^2$ on $\RR$
 it is equivalent to given a differential $1$-form on $\RR.$
 If $f$ is a real-valued $C^2$ map defined on an interval of $\RR$
 on which the derivative vanishes nowhere, we may define the
 following differential $1$-form:

 $$
 \eta(f)=\frac{f''}{f'}dt
 $$

 It follows from the previous definition that
 $$
 \eta(f\circ g)=g^{\ast}\eta(f)+\eta(g)
 $$
 where $g^\ast\eta(f)$ is the pull-back of $\eta(f)$ by $g$:

$$g^{\ast}\eta(f) = \frac{f'' \circ g}{f'\circ g} g' dt$$

 On the other hand, the maps $f$ satisfying $\eta(f)=0$ are
 characterized as the restrictions of affine maps
 (that is, of the form $t\to \lambda t+b$ ).
 Consequently, if $g$ is affine, we have
 $$
 \eta(g\circ f)=\eta(f).
 $$
 Hence, there is a correspondence between a differential $1$-form
 on an interval of $\RR$ and an affine structure on this
 interval. In fact, if $\omega(t)dt$ is a differential $1$-form on an
 interval $I$, then the differential equation $\eta(f)=\omega dt$
 has local solutions which are local diffeomorphisms between $I$
 and an open set of $\RR $.
 Moreover, two of these diffeomorphisms differ by right composition
 by an affine diffeomorphism. Then, the family of this local solutions is a
 system of affine charts on $I.$

Conversely, if $(U_i,f_i)_i$ is a system of charts that defines
 an affine structure of class $C^2,$ then the differential $1$-form
 defined by $\omega(t)dt=\eta(f_i)$ if $t\in U_i$ is independent of the
 choice of $U_i\ni t.$

 \vspace{.5cm}
 We consider a $C^\infty$ Anosov action of $\RR^k$
 on $M$ whose stable foliation $\cF^s$ is of codimension one.
 Then, each leaf of $\cF^{uu}$ is $C^\infty$ diffeomorphic to $\RR.$
 We may assume that $\cF^{uu}$ is orientable, otherwise we consider
 the double covering of $M$. Consequently, it is possible to
 parametrize $\cF^{uu}$ by $u:\RR \times M\to M,$ an application
 such that the signed distance of $u(t,x) \in \cF^{uu}(x)$ to
 $x$ is $t.$ Here we consider the induced metric on $\cF^{uu}.$

 \begin{lemma}
 \label{lem:para unst}
 The application $u$ is continuous. For $x$ fixed, the application
 $u_x:\RR \to M$ defined by $u_x(t)=u(t,x)$ is $C^\infty.$
 Furthermore, the derivatives $\frac{\partial}{\partial t^{\ell}}u(0,x),\ \ell \in \NN,$
 depend continuously on $x.$
 \end{lemma}
 \begin{proof}
 Since $u$ is a flow, it is sufficient to establish the lemma for
 small values of $t.$

 Let $C^{\infty}(\RR,M)$ be the space of
 $C^\infty$ immersions of $\RR$ in $M$ provided of the $C^\infty$
 uniform topology.
 It follows from the theory developed in \cite{hirpush} that for all $x\in M,$
 there exist an open neighborhood $U$ of $x$ and a continuous application
 $\mathcal{I}:U\to C^{\infty}(\RR,M),$ such that, for every $y\in U$,
 the immersion $\mathcal{I}_y=\mathcal{I}(y)$ is a diffeomorphism
 between $\RR$ and a neighborhood of $y$ in $\cF^{uu}$.
 There exist $t_0>0$ such that, for all $(t,y)\in (-t_0,t_0)\times U$ we have that
 $u(t,y)\in \mathcal{I}_y(\RR)$. Hence
 $u(t,y)=\mathcal{I}_y(s(t,y))$ where $s(t,y)\in \RR$ is defined by
 equation
 $$
 \int_0^{s(t,y)}\|\mathcal{I}'_y(\alpha)\|d\alpha=t
 $$
 As $\|\mathcal{I}'_y(\alpha)\|$ is continuous with respect to $y$
 and smooth with respect to $\alpha$, it follows that, $s(t,y)$
 is continuous with respect to $y$ and smooth with respect to $t.$
 This proves that $u$ is continuous and $u_x$ is of class $C^\infty.$
 The last statement of the lemma is trivial.
 \end{proof}

 For a continuous application $\omega:M \to \RR$, the parametrization $u$
 of $\cF^{uu}$ permits us to associate affine structures on the
 leaves $\cF^{uu}(x)$ which are defined by the differential $1$-form
 $\omega(u(t,x))$. This structure will be called of \textit{affine structure along the leaves
 of} $\cF^{uu}$ defined by $\omega$.

 We say that an affine structure along the leaves of $\cF^{uu}$ is
 \textit{invariant by the action} $\phi$ if,  for each $v\in \RR^k$, the
 application $\phi^v|_{\cF^{uu}(x)}:\cF^{uu}(x)\to \cF^{uu}(\phi^v(x)),$
 $x\in M,$ is an affine diffeomorphism.

 \begin{theorem}
 \label{thm:AfinStruct}
 Let $\phi$ be a codimension one Anosov action on $M$ and suppose that
 $\cF^{uu}$ is one dimensional.
 There exists an unique affine structure along the leaves of $\cF^{uu}$
 depending continuously on the points and invariant by the action $\phi.$
 \end{theorem}

 \begin{proof}
 For each $(v,x)\in \RR^k\times M$, let $\tau^v_x:\RR\to \RR$ be the
 application defined by:
 $$
 \phi^v(u(t,x))=u(\tau^v_x(t),\phi^v(x))
 $$
 We claim that a continuous application $\omega:M\to \RR$ defines an invariant affine
 structure along the leaves of $\cF^{uu}$ if and only if
 $$
 \omega=n^v+\delta^v\omega\circ \phi^v=\mathcal{A}^v(\omega),\ \textrm{for all}\ v\in
 \RR^k,
 $$
 where $\delta^v\omega(x)=(\tau^v_x)'(0)$ and $n^v(x)=(\tau^v_x)''(0)/(\tau^v_x)'(0).$
 In fact, $f_y,\ y=\phi^v(x)$ is an affine chart on $\cF^{uu}(\phi^v(x))$
 if only $f_y\circ\tau^v_x$ is an affine chart on $\cF^{uu}(x)$.
 Equivalently:
 $$
 \omega(u(\tau^v_x(t),y))=\frac{f_y''(\tau^v_x(t))}{f_y'(\tau^v_x(t))}=
 \omega\circ\phi^v(u(t,x))\ \Leftrightarrow \
 \omega(u(t,x))=\frac{(f_y\circ\tau^v_x)''(t)}{(f_y\circ\tau^v_x)'(t)}
 .
 $$
 Hence, as $\tau^v_{u(t,x)}(s)=\tau^v_x(t+s)$, we obtain:
 $$
 \begin{array}{cll}
 \omega(u(t,x))&= &{\displaystyle \frac{(f_y\circ\tau^v_x)''(t)}{(f_y\circ\tau^v_x)'(t)}}\\
 & = &{\displaystyle \frac{(\tau^v_x)''(t)}{(\tau^v_x)'(t)}+(\tau^v_x)'(t)\frac{f_y''(\tau^v_x(t))}{f_y'(\tau^v_x(t))}}\\
    & = &{\displaystyle \frac{(\tau^v_{u(t,x)})''(0)}{(\tau^v_{u(t,x)})'(0)}+(\tau^v_{u(t,x)})'(0)\omega\circ\phi^v(u(t,x))} \\
    & = & n^v(u(t,x))+\delta^v(u(t,x))\omega\circ\phi^v(u(t,x))\\
    & = & \mathcal{A}^v(\omega)(u(t,x))
 \end{array}
 $$
 This proves our claim.

 The applications $\mathcal{A}^v$ acting on the Banach space of the
 continuous applications of $M$ on $\RR$ provided of the uniform
 norm. By definition of Anosov action, if $v$ is an element of the
 Anosov chamber, we have that $\delta^{sv}, \ s<0$ has uniform norm less
 that one. This implies that $\mathcal{A}^{sv},\ s<0$ is a contraction, hence
 $\mathcal{A}^{sv},\ s<0$ admit an unique fixed point. Finally,
 since $\mathcal{A}^v \circ \mathcal{A}^w=\mathcal{A}^w\circ
 \mathcal{A}^v$ for all $v,w\in \RR^k$, there exists an unique fixed point $\omega$
 for all $\mathcal{A}^v$. This finishes the proof.
  \end{proof}

 Real affine structures on the real line are well-known: they are all affinely isomorphic to
 the segment $(0, 1)$, the half-line $(0,+\infty)$, or the complete affine line $(-\infty, +\infty)$.
 In the latter case, the affine structure is said {\em complete.}

 \begin{lemma}
 \label{le:complete}
 Every leaf of $\cF^{uu}$, endowed with the affine structure provided by Theorem~\ref{thm:AfinStruct},
 is complete.
 \end{lemma}

 \begin{proof}
 For every $x$ in $M$, there is a unique affine map $f_x: \cF^{uu}(x) \to \RR$ mapping
 $x$ on $0$ and the point $u(1,x)$ at distance $1$ on $1$. The image of $f_x$ is an
 interval $(\alpha(x), \beta(x))$. We aim to prove that $\alpha(x)=-\infty$ and $\beta(x)=+\infty$.

 For every $v$ in $\RR^k$ and $x$ in $M$, the restriction of $\phi^v$ on
 $\cF^{uu}(x)$ induces an affine transformation, even linear, of the affine line, of the form
 $z \to \lambda(v,x)z$:

 $$f_{\phi^v(x)} \circ \phi^v = \lambda(v,x)f_x$$

 Hence:

 $$\beta(\phi^v(x))=\lambda(v,x)\beta(x)$$

 Recall that:

 $$
 \phi^v(u(t,x))=u(\tau^v_x(t),\phi^v(x))
 $$

 Since by definition $f_x(u(1,x))=1$, we get:

 $$
 \lambda(v,x)=f_{\phi^v(x)}(u(\tau^v_x(t),\phi^v(x)))
 $$

 Hence, for $v=ta$, where $a$ is the codimension Anosov element, $\lambda(-ta,x)$ is arbitrarily
 small if $t>0$ is sufficiently big. Therefore, if $\beta(x)$ is not $+\infty$, $\beta(\phi^{-ta}(x))$
 takes arbitrarily small value. This is a contradiction since obviously $\beta>1$ everywhere.
 Therefore, $\beta$ is infinite everywhere.

 The proof of $\alpha=-\infty$ is similar.
 \end{proof}

 %We begin by remembering that there exist a correspondence between projective
% structure on $\RR$ and quadratic differential on $\RR.$
% If $f$ is a real-valued $C^3$ map defined on an interval of $\RR$
% on which the derivative is never annulled, we may define the
% Schwarzian derivative of $f$ as:
%
% $$
% S(f)=\frac{f'''}{f'}-\frac{3}{2}\Big(\frac{f''}{f'}\Big)^2
% $$
% We introduce the quadratic differential $s(f)$ defined by:
% $$
% s(f)=S(f)dt^2
% $$
% It follows from the previous definition that
% $$
% s(f\circ g)=g^{\ast}s(f)+s(g)
% $$
% where $g^{\ast}s(f)$ is the pullback by $g$ of $s(f)$, that is
% $$
% g^{\ast}s(f)=(g')^2S(f)\circ g dt^2.
% $$
% On the other hand, the maps $f$ satisfying $s(f)=0$ are
% characterized as the restrictions of projective diffeomorphisms
% (that is, of the form $t\to \frac{at+b}{ct+d}, \ ad-bc\neq 0$ ).
% Consequently, if $g$ is projective, we have
% $$
% s(g\circ f)=s(f).
% $$
% Hence, there is a correspondence between a quadratic differential
% on an interval of $\RR$ and a projective structure on this
% interval. In fact, if $a(t)dt^2$ is a quadratic differential on an
% interval $I$, then the differential equation $S(f)=a$
% has local solutions which are local diffeomorphisms between $I$
% and an open set of $P_1=\RR \cup \{\infty\}$.
% Moreover, two of these diffeomorphisms differ by composing , in the sink,
% by a projective map. Then, the family of this local solutions is a
% system of projective charts on $I.$

 \subsection{Irreducible codimension one Anosov actions}
 A codimension one Anosov action $\phi$ of $\RR^k$ on $M$ is said to be
 \textit{irreducible} if for any $v\in \RR^k-\{0\}$ and $x\in M$
 with $\phi^v(x)=x$ we have that ${\rm Hol}_{\gamma}$, the holonomy along
 of $\gamma=\{\phi^{sv}(x);\ s\in [0,1]\}$ of $\cF^{s}(x)$, is
 a topological contraction or a topological expansion.

 \begin{remark}
 \label{rk:defreduc}
 {\rm
 It follows from Theorem~\ref{thm:AfinStruct} that the holonomy along $\gamma$ is differentially linearizable. Therefore,
 an equivalent definition of irreducibility is to require that the holonomy along
 $\gamma$ is non-trivial.
 }
 \end{remark}

 \begin{remark}
 When $k=1$, the case that the action is a flow, all the
 codimension one Anosov actions are irreducibles.
 \end{remark}

 \begin{theorem}
 \label{thm:ReducingAction}
 Let $\phi:\RR^k\times M\to M$ be a codimension one Anosov action.
 Then, there exists a free abelian subgroup $H_0 \approx \RR^\ell$ of  $\RR^k$, a lattice
 $\Gamma_0 \subset H_0$, a smooth
 $(n-\ell)$-manifold $\bar{M}$, and $p:M\to \bar{M}$ a smooth
 $\TT^{\ell}$-principal bundle such that:
 \begin{enumerate}
   \item $\Gamma_0$ is the kernel of $\phi;$
   \item every orbit of ${\phi}_0=\phi|_{H_0 \times M}$ is a fiber of $p:M\to \bar{M}$.
   In particular, $\bar{M}$ is the orbit space of $\phi_0$;
   \item $\phi$ induces an irreducible codimension one Anosov action
   $\bar{\phi}: \bar{H} \times \bar{M} \to \bar{M}$ where $\bar{H}=\RR^k/H_0$.
    \end{enumerate}
 \end{theorem}

 The proof of Theorem~\ref{thm:ReducingAction} essentially relies on the following lemma:

 \begin{lemma}
 \label{le:trivial}
 Let $v$ be an element of $\RR^k$ and $x$ an element of $M$ such that $\phi^v(x)=x$.
Then, either $x$ is a repelling of attracting (and therefore, unique) fixed point of the restriction of
$\phi^v$ to $\cF^{uu}(x)$, or the action of $\phi^v$ on the entire manifold $M$ is trivial.
 \end{lemma}

 \begin{proof}
 For every $x$ in $M$ and every $w$ in $\RR^k$ such that $\cF^{ss}(x)=\phi^w(\cF^{ss}(x))$ we consider any loop
 in $\cF^{s}(x)$ which is the composition of  $t\in [0,1] \to \phi^{tw}(y)$ with any path in $\cF^{ss}(x)$ joining
 $\phi^w(x)$ to $x$. Since $\cF^{ss}(x)$ is a plane (Remark~\ref{rk.zero}) all these loops are homotopic one to the other
 in $\cF^s(x)$; in particular, the  holonomy of $\cF^s$ along any of them is well-defined and does not depend on the loop.
 We denote it by $h^w_x$.

 Assume that  $x$ is the one appearing in the statement of the lemma.
 Then according to Theorem~\ref{thm:AfinStruct}, the restriction of $\phi^v$ to $\cF^{uu}(x)$ is conjugated to
 an affine transformation of the real affine line. Therefore, in order to prove the theorem, we just have to consider
 the case where $h^v_x$ is trivial.

We define $\Omega_w$ as the set comprising the points $x$ in $M$ such that
$\phi^w(x) \in \cF^{ss}(x)$ and for which the holonomy $h^w_x$ is trivial. This set is obviously
$\phi$-invariant. By the discussion above,
we can assume that $\Omega_v$ is non empty. Moreover, for every $x$ in $\Omega_w$, and every $y$ in $\cF^{ss}(x)$, the loops
considered above associated to respectively $x$, $y$,
are freely homotopic one to the other in $\cF^s(x)$. Hence $h^w_y=h^w_x$. It follows
that $\Omega_w$ is saturated by $\cF^s$.

Finally, for every subset $U$ of $\RR^k$, let
$\Omega_U$ be the union of the $\Omega_w$ for $w$ in $U$.
For every $x$ in $\Omega_v$, since the holonomy $h^v_x$ is trivial, for every $y$ in $\cF^{uu}(x)$
near $x$ the point $\phi^v(y)$ lies in the local stable leaf of $y$. It follows that $y$ lies in $\Omega_w$
for some $w$ close to $v$ in $\RR^k$. Since $M$ is compact, for every neighborhood $U$ of $v$ in
$\RR^k$, there exists $\delta>0$ such that every $y$ in $M$, lying on a local unstable leaf
$\cF^{uu}_\delta(x)$ with $x$ in $\Omega_v$, belongs to $\Omega_U$.

Now, at the one hand we know that $\Omega_U$ is $\phi$-invariant. On the other hand, since the $\cF^{uu}$-saturation
of any $\cF^{s}$-invariant subset is the entire $M$, for every $y$ in $M$ the point $\phi^{ta}(y)$ lies in $\cF^{uu}_\delta(x)$
for some $t<0$, where $x$ is an element of $\Omega_v$. It follows that $\Omega_U$ is the entire $M$. Since $U$
is arbitrary, we get the equality $M=\Omega_v$.

Consider now a compact $\phi$-orbit $\cO$.  For some $\delta>0$ and every $x$ in
$\cO$, the intersection $\cO \cap \cF^{ss}_\delta(x)$ is reduced to $x$. For every $y$ in
$\cF^{ss}(x)$, and for every $t>0$ sufficiently big, $\phi^{ta}(y)$ belongs to $ \cF^{ss}_\delta(\phi^{ta}(x))$.
Hence $\cF^{ss}(x) \cap \cO = \{ x \}$. It follows that every point in ${\rm Comp}(\phi) \cap \Omega_v = {\rm Comp}(\phi)$
is fixed by $\phi_v$. Hence, the restriction of $\phi^v$ to the closure $\Omega(\phi)$ of ${\rm Comp}(\phi)$ is trivial.

Finally, assume that $x$ is an arbitrary element of $M=\Omega_v$. Let $(t_n)_{(n \in \NN)}$ be a sequence
of positive real number diverging to $+\infty$ and such that $x_n=\phi^{-t_na}(x)$ converges to some element $x_\infty$
of $\Omega(\phi) \subset {\rm Fix}(\phi)$. Then, due to the proximity to $x_\infty$,
for every $\epsilon>0$, and for $n$ is sufficiently big, there is a path
$c_n$ in $\cF^{ss}(x_n)$ of length $\leq \epsilon$ joining $x_n$ to $\phi^v(x_n)$.
Then, $\phi^{t_na}(c_n)$ is a path of length $\leq \epsilon$ joining $x$ to $\phi^v(x)$. Since $\epsilon$
is arbitrary, we get $\phi^v(x)=x$. This achieves the proof of the lemma.
  \end{proof}

{\emph{Proof of Theorem~\ref{thm:ReducingAction}. }
 Let $\Gamma_0$ be the kernel of $\phi$.
 Since the action is locally free,  $\Gamma_0$ is discrete, isomorphic to
 $\mathbb{Z}^{\ell}$ for some integer $\ell \geq 0$. Let $H_0$ be the subspace of $\RR^k$
 generated by $\Gamma_0$. Observe that for
 every $x$ in $M$, and every $v$ in $\Gamma_0$, the holonomy $h^v_x$ is well-defined, and trivial
 (cf. the notations introduced in the proof of Lemma~\ref{le:trivial}).

 The torus $\TT^\ell= \Gamma_0\backslash{H_0}$ acts on $M$; since it is compact, this action is proper.
Moreover, this action is free: indeed, if some $v$ in $H_0$ fixes
some $x$, then the holonomy $h^v_x$ is trivial since $v$ is a linear
combination of elements of $\Gamma_0$ (or, better to say, since some
of the iterates $nv$, for integers $n$, are arbitrarily approximated
by elements of $\Gamma_0$). According to Lemma~\ref{le:trivial}, the
action of $v$ on $M$ is trivial, i.e. $v$ belongs to $\Gamma_0$.

Therefore, the quotient space $\bar{M}$ is a closed $(n-\ell)$-dimensional manifold, and
the quotient map $p: M \to \bar{M}$ is a principal $\TT^\ell$-bundle. The action
of $\RR^k$ on $M$ induces an action of $\bar{H}= \RR^k/H_0$. It is straightforward to check that this action is
Anosov, and of codimension one.

Finally, if $\bar{\phi}$ is not irreducible, there is a non-trivial element $\bar{v}$ of $\bar{H}$ fixing a point $\bar{x}$ in
$\bar{M}$ and such that $h^{\bar{v}}_{\bar{x}}$ is  trivial. There is a representant $v$ of $\bar{v}$ in $\RR^k$ fixing a  point $x$ in $M$
above $\bar{x}$, and such that $h^{v}_x$ is trivial. According to
Lemma~\ref{le:trivial}, $v$ belongs to $\Gamma_0 \subset H_0$. Hence $\bar{v}$ is trivial. This contradiction achieves the proof of
the theorem.
\fin

 \begin{remark}
 \label{rem:irreducible}
 Let $\phi:\RR^k\times M\to M$ be a codimension one Anosov action and
 $\bar{\phi}$ an action of $\bar{H}$ (that is isomorphic to $\RR^{k-\ell}$) on $\bar{M}$ as
 in Theorem \ref{thm:ReducingAction}. Then:
 \begin{enumerate}
   \item $\dim M > k+2 \Longleftrightarrow \dim \bar{M} > \dim \bar{H} +2$;
   \item $\phi$ is transitive $\Longleftrightarrow \bar{\phi}$ is
   transitive;
   \item $\phi$ is irreducible $\Longleftrightarrow \bar{H}=\RR^k.$
 \end{enumerate}
 \end{remark}

 \subsection{The orbit space of an irreducible codimension one Anosov action}

 Let $\pi:\widetilde{M}\to M$ be the universal covering map of $M$ and $\widetilde{\phi}$ be the
 lift of $\phi$ on $\widetilde{M}$. The foliations $\cF^{ss},\ \cF^{uu},\ \cF^{s}$ and
 $\cF^{u}$ lift to foliations $\widetilde{\cF}^{ss},\ \widetilde{\cF}^{uu},\ \widetilde{\cF}^{s}$ and
 $\widetilde{\cF}^{u}$ in $\widetilde{M}$. We denote by $Q^{\phi}$ be the
 orbit space of $\tilde{\phi}$ and $\pi^{\phi}:\widetilde{M}\to Q^{\phi}$ be the
 canonical projection. This section is devoted to the proof of the following theorem, which is
 a keystone of the proof of the main theorem.

 \begin{theorem}
 \label{thm:orbitspace}
 If $\phi$ is an irreducible codimension one Anosov action of $\RR^k$ on $M$ then
 $Q_{\phi}$, the orbit space of $\widetilde{\phi}$, is homeomorphic to $\RR^{n-k}.$
 \end{theorem}

 Let $p=n-1-k$. Since by hypothesis $\cF^{uu}$ is one dimensional,
 $p$ is the dimension of $\cF^{ss}$.

 \begin{proposition}
 \label{pro.transverse}
 Every loop in $M$ transverse to $\cF^s$ is homotopically non-trivial in $M$.
  \end{proposition}

 \begin{proof}
 By a Theorem of Haefliger (see \cite[Proposition 7.3.2]{conlon}), if some transverse loop is
 homotopically trivial, then there is a leaf ${F}$,
 of  ${\cF}^{s}$ containing a loop $c: [0, 1] \to  {F}$, homotopically non-trivial in ${F}$, and
 such that the holonomy of  ${\cF}^{s}$ along $c$ is trivial on one side.

 According to Remark~\ref{rk.zero} there is a fibration $p_F: F \to P$ where $P$ is a flat
 cylinder, and whose fibers are strong stable leaves. Moreover, the restriction of
 $p_F$ to the orbit $\cO_x$ of $x=c(0)$ is a covering map. Hence we can lift in $\cO_x$
 the curve $p \circ c$. In other words, there is a continuous path
 $c': [0, 1] \to \RR^k$ such that for every $t$ in $[0, 1]$ the
 image $c'(t)$ lies in the same leaf of  ${\cF}^{ss}$ than $\phi^{c(t)}(x)$. Let $v=c'(1)$: since $c$
 is homotopic to the loop obtained by composing $c'$ with any path in $\cF^{ss}(c(0))$ joining
 $c'(1)$ to $c(1)=c(0)=x$, we get that the holonomy $h^v_{x}$ is trivial on one side. Since
 it is linearizable, $h^v_{x}$ is trivial. According to Lemma~\ref{le:trivial}, and since
 $\phi$ is irreducible, it means that $c'(1)=0$. Hence $c'$ is homotopically trivial in
 $\RR^k$. Hence $p_F \circ c'$ is homotopically trivial in $P$. But since $p_F$ is a trivial
 fibration with contractible leaves and $p_F \circ c = p_F \circ c'$, we obtain that $c$ is
 homotopically trivial in $F$. Contradiction.
 \end{proof}

 \begin{corollary}
 \label{cor:incompressible}
 The orbits of ${\phi}$ are incompressible: every loop in a $\phi$-orbit $\cO$ which
 is homotopically non-trivial in $\cO$ is homotopically non-trivial in $M$.
 \end{corollary}

\begin{proof}
Every loop in $\cO$ is homotopic to a trajectory $t \to c(t)=\phi^{tv}(x); \  t\in [0, 1], v \in \RR^n, \phi^v(x)=x.$
Since $\phi$ is irreducible, the holonomy of $\cF^s$ along $c$ is non-trivial. It follows
that there is a loop homotopic to $c$ and transverse to $\cF^s$. The corollary follows from
Proposition~\ref{pro.transverse}.
\end{proof}

 A foliation is said to be \textit{by closed planes} if all the leaves are closed and images of embeddings of
$\RR^n.$

 \begin{corollary}
 \label{cor:closedleaves}
 Let $\phi$ be an irreducible codimension one Anosov action on $M$. The
 foliations $\widetilde{\cF}^{uu}$, $\widetilde{\cF}^{ss}$, $\widetilde{\cF}^{u}$,
 $\widetilde{\cF}^{s}$ and the foliation defined by $\widetilde{\phi}$ are by closed planes.
  The intersection between a leaf of $\widetilde{\cF}^{u}$ and a leaf of
 $\widetilde{\cF}^{s}$ is at most an orbit of $\widetilde{\phi}$.
 Every orbit of $\widetilde{\phi}$ meets a leaf of $\widetilde{\cF}^{uu}$ or $\widetilde{\cF}^{ss}$
 at most once.
 \end{corollary}

 \begin{proof}
 According to Corollary~\ref{cor:incompressible}, $\tilde{\phi}$ is a free action. Moreover,
 if $\widetilde{\cF}^{ss}(x) = \widetilde{\cF}^{ss}(\tilde{\phi}^v(x))$ for some non-trivial $v$,
 then $h^v_{\pi(x)}$ is non-trivial (since $\phi$ is irreducible), but is also
 the holonomy of $\cF^s$ along a closed loop in $\cF^s(\pi(x))$ homotopically trivial in $M$.
 It is in contradiction with Proposition~\ref{pro.transverse}.

 Therefore, every orbit of $\tilde{\phi}$ intersects every leaf of $\widetilde{\cF}^{ss}$ at most once.
 Since every leaf of $\widetilde{\cF}^s$ is the saturation under $\tilde{\phi}$ of a leaf of
 $\widetilde{\cF}^{ss}$, it follows that it is an injective immersion of $\RR^{p+k}$. Similarly,
 every leaf of $\widetilde{\cF}^u$ is an injective immersion of $\RR^{2+k}$.

 It is easy to show, by a standard argument, that if a leaf $\widetilde{\cF}^s$ is not closed, then
 there is a loop in $\widetilde{M}$ transverse to $\widetilde{\cF}^s$, giving a contradiction with
 Proposition~\ref{pro.transverse}. Hence $\widetilde{\cF}^s$ is a foliation by closed planes.
 The statement for leaves of $\widetilde{\cF}^{ss}$ follows.

 Since $\widetilde{M}$ is simply connected, every foliation in it are oriented and transversely oriented.
 Being closed hypersurfaces, leaves of $\widetilde{\cF}^s$ disconnect $\widetilde{M}$. Hence, a leaf
 of $\widetilde{\cF}^{uu}$ intersecting a leaf $F$ of $\widetilde{\cF}^{s}$ enters in one side of $F$,
 and cannot cross $F$ once more afterwards, due to orientation considerations. In other words, every
 leaf of $\widetilde{\cF}^{uu}$ intersects every leaf of $\widetilde{\cF}^{s}$ at most once. In particular,
 it cannot accumulate somewhere, i.e. it is closed.

 In order to achieve the proof of the corollary, we just have to prove that the $\tilde{\phi}$-orbits are
 closed. But this is clear, since each of them is the intersection between a weak stable leaf and a
 weak unstable leaf, that we have shown to be closed.
 \end{proof}

 Consequently, by a theorem of Palmeira \cite{palmeira}:

 \begin{corollary}
 With the same hypotheses of above proposition, the universal
 covering of $M$ is diffeomorphic to $\RR^n.$
 \fin
 \end{corollary}

 \begin{lemma}
 If $\phi$ be a codimension one Anosov action on $M$, then the orbit
 space of $\widetilde{\phi}$ is Hausdorff.
 \end{lemma}
 \begin{proof}
 By contradiction, we assume that there exist two different $\widetilde{\phi}$-orbits
 $\widetilde{\mathcal{O}}_{x_1}$ and $\widetilde{\mathcal{O}}_{x_2}$
 which are non-separable. Then, the saturation by $\widetilde{\cF}^{uu}$ of $\widetilde{\cF}^s(x_1)$
 and $\widetilde{\cF}^s(x_2)$ are two non disjoint neighborhoods of
 $\widetilde{\mathcal{O}}_{x_1}$ and $\widetilde{\mathcal{O}}_{x_2}$, respectively.

First, we assume that $\widetilde{\mathcal{O}}_{x_1}$ and
$\widetilde{\mathcal{O}}_{x_2}$ are contained in the same leaf of
$\widetilde{\cF}^s$. Hence, we can assume that
$\widetilde{\cF}^{ss}(x_1)=\widetilde{\cF}^{ss}(x_2)=F_0$. Let $U_1$
and $U_2$ be the disjoint neighborhoods in $F_0$ of $x_1$ and $x_2,$
respectively. It follows, from Corollary \ref{cor:closedleaves},
that the saturation by $\widetilde{\cF}^{u}$ of $U_1$ and $U_2$ are
two  disjoint $\widetilde{\phi}$-invariant neighborhoods of
$\widetilde{\mathcal{O}}_{x_1}$ and $\widetilde{\mathcal{O}}_{x_2}$.
This contradicts our assumption.

Hence, $\widetilde{\cF}^s(x_1)\neq \widetilde{\cF}^s(x_2).$ The saturation
by $\widetilde{\cF}^{uu}$ of $\widetilde{\cF}^s(x_1)$ and $\widetilde{\cF}^s(x_2)$
cannot be disjoint since they are neighborhoods of respectively $x_1$, $x_2.$
There exist
 $y_1\in \widetilde{\cF}^s(x_1)$ and $y_2\in \widetilde{\cF}^s(x_2)$
 such that $\widetilde{\cF}^{uu}(y_1)$=$\widetilde{\cF}^{uu}(y_2)$.
 Since $y_1\neq y_2$, there exist disjoint
 neighborhoods $U_1$ and $U_2$  in $\widetilde{\cF}^{uu}(y_1)$ of $y_1$ and $y_2$,
 respectively. The saturation by $\widetilde{\cF}^{s}$ of $U_1$ and $U_2$ are
 two $\widetilde{\phi}$-invariant neighborhoods of $\widetilde{\mathcal{O}}_{x_1}$
 and $\widetilde{\mathcal{O}}_{x_2}$ which, by our assumption, are
 non disjoint. In this case, a leaf of $\widetilde{\cF}^{s}$ passing
 by a point in the intersection of these neighborhoods meet
  $\widetilde{\cF}^{uu}(y_1)$ in two points, this contradicts the
  Corollary \ref{cor:closedleaves} and finishes the proof.
 \end{proof}

 \vspace{.5cm}
 \begin{proof}[Proof of Theorem \ref{thm:orbitspace}]
 Let $x\in \widetilde{M}$ and $U$ be a neighborhood of $x$
 which is a lift of a product neighborhood in the sense of Theorem
 \ref{thm:local product}. Let $\Sigma \subset U$ be a smooth
 $(n-k)$-submanifold which is transverse to $\phi$. First we are going
 to show that every orbit of $\phi$ meets $\Sigma$ in at most one point.
 By contradiction, we assume that there exists an orbit $\cO$
 that meets $\Sigma$ in two points. As $\cO$ meets a leaf of
 $\widetilde{\cF}^{ss}$ in at most one point (Corollary \ref{cor:closedleaves}),
 then $\cO$ meets $U$ along two different leaves of
 $\widetilde{\cF}^s|_U$. Hence, since $U$ is a product neighborhood,
 we have that there are leaves of $\widetilde{\cF}^{uu}$ that meets
 a leaf of $\widetilde{\cF}^{ss}$ at two points. But this is impossible
 by Corollary \ref{cor:closedleaves}.

 Let $(\Sigma_i)_{i\in I}$ be a family of transverses as above whose
 union meets all the orbits of $\phi.$ Then $\{(\Sigma_i,\pi^{\phi}|_{\Sigma_i}); \ i\in I\}$
 defines a differentiable structure on $Q^{\phi}$ whose class of
 differentiability is the same that of the action $\phi.$
 Moreover, $\pi^{\phi}$ is a  locally trivial bundle. Thus $Q^{\phi}$
 is a manifold of dimension $n-k$, Hausdorff and, as shown by the
 exact sequence of homotopy groups for the bundle $\pi^{\phi}:\widetilde{M}\to
 Q^{\phi}$, simply connected.  Since $\widetilde{\cF}^s$ induces on
 $Q^{\phi}$ a codimension one foliation by planes, we conclude, once more
 by Palmeira's Theorem, that
 $Q^{\phi}$ is diffeomorphic to $\RR^{n-k}$.
 \end{proof}

 \begin{proposition}
 \label{pro:OrbIrrecAct}
 All the non compact orbits of an irreducible codimension one Anosov action are
 planes.
 \end{proposition}
 \begin{proof}
 Let $\cO$ be an orbit of an irreducible codimension one Anosov
 action $\phi:\RR^k \times M\to M$. Suppose that $\cO$ is not
 a plane, i.e. that  $\phi^v(x)=x$ on $\cO$ for some  $v\in \RR^k -\{0\}$. Let $y\in \overline{\cO}$.
 Then there exist a sequence $\{x_n\}$ of elements of $\cO$
 such that $x_n \to y$. Thus,
 $\phi^v(y)=y$.  Since the action is irreducible, the holonomy $h^v_y$ is non trivial. It follows that
all the $x_n$, for $n$ sufficiently big, lie in the same local stable leaf $\cF^{s}_{\delta}(y)$.

Hence the closure $\overline{\cO}$ in $M$, which is compact, is contained in the weak stable leaf $F=\cF^{s}(y)=\cF^{s}(x).$
It follows that the space $P=p_F(\cO)=p_F(\overline{\cO})$ of strong stable leaves in $F$ is compact (cf. the notations in Remark~\ref{rk.zero}),
hence $\Gamma_L$ is a lattice in $\RR^k$. According to Lemma~\ref{le:compactcausal}, $\Gamma_L$
contains an Anosov element $b$ contained in the chamber $\cA_a$.
The restriction of $\phi^a$ to the strong stable leaf $L$ is a contraction, hence contains a fixed point $z$.
By Remark~\ref{rk.first}, the $\phi$-orbit of $z$ is compact, and $\Gamma_L$ is the isotropy
group of $z$. Now, since $\overline{\cO}$  is compact, the same is true for the intersection $L \cap \overline{\cO}$.
On the other hand, the negative iterates $\phi^{-na}(z')$ of a point $z'$ in $L$ different from $z$ escape from
any compact subset of $L$. Therefore, $\overline{\cO}$ is the compact orbit of $z$.
 \end{proof}

As an immediate corollary of Theorem~\ref{thm:ReducingAction} and Proposition~\ref{pro:OrbIrrecAct} we get:

\begin{corollary}
 Let $\phi:\RR^k\times M\to M$ be a codimension one Anosov action, not necessarily irreducible.
 Every non-compact $\phi$-orbit is diffeomorphic to $\TT^\ell\times\RR^{k-\ell}$.
\fin
\end{corollary}

 \section{Codimension one Anosov $\RR^k$-actions are transitive}
 \label{sec:proof}
 In this section we prove the main Theorem:

 \vspace{.5cm}
 \noindent
 \textbf{Theorem 2.}
 \textit{Let $\phi$ be a codimension one Anosov action of $\RR^k$ on a closed manifold $M$
 of dimension greater than $k+2$. Then, any regular subcone $\cC$ admits a dense orbit in $M$.}

 \vspace{.5cm}

  In what follows, we consider a codimension one Anosov action
 $\phi:\RR^k\times M\to M$, with $\dim M > k+2$ and a regular
 subcone $\cC$. Furthermore, we will consider
 the foliations $\cF^{ss},\ \cF^{uu},\ \cF^{s}$ and $\cF^{u}$ which
 corresponds to the chamber $\mathcal{A}_a$ containing $\cC.$

 \subsection{Bi-homoclinic points and transitivity}
 Let $u$ be the parametrization of $\cF^{uu}$ which was studied in
 Lemma \ref{lem:para unst}. We define
 $$
 \begin{array}{c}
 \mathcal{H}^+=\{x\in M;\ (x,+\infty)\cap \cF^s(x)=\emptyset\}\\

 \mathcal{H}^-=\{x\in M;\ (-\infty,x)\cap \cF^s(x)=\emptyset\}
 \end{array}
 $$

 where $(x,+\infty)=\{u(\tau,x);\ \tau >0\}$ and $(-\infty,x)=\{u(\tau,x);\ \tau
 <0\}$.

 \begin{definition}
 {\rm
 A point $x\in M$ is said to be \textit{bi-homoclinic} if $x\notin \mathcal{H}^+\cup \mathcal{H}^-.$
 }
 \end{definition}

 The following results establish a criterion for transitivity.

 \begin{proposition}
 \label{prop:pointbihomoclinic}
 If every point in $M$ is bi-homoclinic, then every leaf of $\cF^s$ is dense.
 \end{proposition}

 \begin{proof}
 For all $\tau_0 \in \RR$ we consider the set
 $\{x\in M;\ (x,u(\tau_0,x))\cap \cF^s(x)\neq \emptyset\}$.
 It follows, from Theorem \ref{thm:local product}, that this set
 is open. Hence, since $M$ is compact, there is $\tau_0 \in \RR$
 such that
 $$
 \textrm{for all}\ x \in M, \ \ \textrm{there exists}\ \tau \in (0,\tau_0)\ \textrm{such that}
 \ u(\tau,x)\in \cF^s(x)
 $$
 Similarly, increasing $\tau_0$, if necessary, we have
  $$
 \textrm{for all}\ x \in M, \ \ \textrm{there exists}\ \tau \in (-\tau_0,0)\ \textrm{such that}
 \ u(\tau,x)\in \cF^s(x).
 $$

 Thus, there exists $\ell >0$ such that if $I$ is an interval contained in a
 leaf of $\cF^{uu}$ whose arc length is grater than $\ell$, then each
 leaf of $\cF^s$ intercepting $I$ contains at least three points of $I$.

 On the other hand, any interval contained in a leaf of $\cF^{uu}$
 admits an iterate by $\phi^a$ ($a$ Anosov element) whose arc length is
 greater than $\ell.$ This implies that any interval of
 $\cF^{uu}(x)$ meets $\cF^{s}(x)$, and hence the closure of
 $\cF^{s}(x)$ is an open set. Therefore $\cF^{s}(x)$ is dense in $M.$
 \end{proof}

 \begin{lemma}
 \label{lem:transitivity}
 If every leaf of $\cF^{s}$ is dense in $M$, then $\cC$ admits a dense orbit.
 \end{lemma}
 \begin{proof}
 Let $a\in \cC$ be an Anosov element and $\phi^{ta}$ the
 corresponding flow. According to Remark~\ref{rk.nonwanderingflow}:
 $$
 \Omega(\phi)=\Omega(\phi^{ta}).
 $$
 By a result of Conley \cite{conley}, there exists $L:M\to \RR$ a
 complete Lyapunov function for the flow $\phi^{ta},$ meaning that
 for $t >0$ we have $L(\phi^{ta}(x)) \leq L(x)$, and the equality holds
 only if $x$ lies in $\Omega(\phi^{ta})$. On the other
 hand, Theorem \ref{thm:espdecomp} and the hypotheses imply that
 the $\Omega(\phi)$ admits only one basic block, in particular,
 that there $\cC$ admits a dense orbit in $\Omega(\phi)$. Moreover,
 for every $x$, $y$ in the basic block $\Omega(\phi)$, the unstable leaf
 of $x$ intersects the stable leaf of $y$, and the stable leaf of $x$ intersects
 the unstable leaf of $y$. From the former we get $L(x) \geq L(y)$, and from the latter,
 $L(y) \geq L(x)$. Hence,
 the restriction of $L$ to $\Omega(\phi)$ is constant, say, vanishes.
 Then, for every $x$ in $M$, the inequalities $L(x)\geq 0$ and $L(x)\leq 0$ hold
 (the former because the $\alpha$-limit set of the $\phi^{ta}$-orbit of
 $x$ is non-empty, the latter
 because the $\omega$-limit set is non-empty). It follows that
 $L$ vanishes everywhere, and $\Omega(\phi)=M$.
 Therefore $\cC$ is transitive.
 \end{proof}

 \subsection{Proof of Main Theorem}
 \begin{lemma}
 \label{lem:nonhomoclinic}
 Assume that $\phi$ is an irreducible codimension one Anosov action.
 The sets $\mathcal{H}^-$ and $\mathcal{H}^+$ are unions of
 compact orbits.
 \end{lemma}
 \begin{proof}
 Note that $\mathcal{H}^-$ and $\mathcal{H}^+$ are closed
 invariant sets. We will show the lemma for $\mathcal{H}^+$, the
 case of $\mathcal{H}^-$ is analogous.

 Cover $M$ by a finite collection $(U_i)_{1 \leq i \leq N)}$ of
 product neighborhoods as in Theorem \ref{thm:local product}.
 We claim that the intersection of any orbit in $\cH^+$ with
 every $U_i$ is connected. This will show the lemma. Indeed, since the $U_i$ are in finite number,
 the orbit under consideration is compact.
 Moreover, since the area of local orbits contained in every $U_i$
 is uniformly bounded from above, it also implies that there is an uniform bound
 on the area of compact orbits in $\cH^+$. According to Lemma~\ref{le:raghu},
 $\mathcal{H}^+$ is the union of a  finite number of compact orbits.

 We are going to show our claim above. Let $y_0\in \mathcal{H}^+$, and
 let $F$ be its weak stable leaf. By the very definition of $\mathcal{H}^+$,
 for every $i$, the intersection between the orbit $\cO_0$ of $y_0$ and the product neighborhood
 $U_i$ must be contained in a single plaque $F_i$. The union of the closure of the $F_i$
 is compact; it follows that $\cO_0$ is relatively compact, not only in the manifold $M$, but also
 in the leaf $F$ equipped with its own leaf topology.

 Recall that there is there is a bundle map $p_F: F \to P$, whose restriction to $\cO_0$
 is a covering map (remark~\ref{rk.zero}). Since $\cO_0$ is relatively compact in $F$,
 the base manifold $P=\RR^k/\Gamma_F$ is compact. By Lemma~\ref{le:compactcausal},
 the lattice $\Gamma_L$ must contain an element $v$ of $\cA_a$. The restriction
 of $\phi^v$ to $\cF^{ss}(y_0)$ is then a contracting map, hence admitting
 a fixed point $y_1$. Moreover, since $\cO_0$ is relatively compact, the leaf distance
 between the iterates $\phi^{-tv}(y_0)$ and $\phi^{-tv}(y_1)$ is bounded.
 It is possible if and only if $y_0=y_1$, in particular, $y_0$ is fixed by the Anosov element
 $\phi^v$.

 According to item (4) of Remark~\ref{rk.first}, the orbit $\cO_0$ is compact.
 The claim and the lemma then follow,
 since $F$ contains at most one compact orbit, and that this compact orbit intersects every
 strong stable leaf at most once.
 \end{proof}

 \vspace{.5cm}
 \begin{proof}[Proof of Main Theorem]
 As by Remark \ref{rem:irreducible} we can assume that $\phi$ is
 irreducible, then by Proposition \ref{prop:pointbihomoclinic} and
 Lemma \ref{lem:transitivity}, it is sufficient to show that the sets
 $\mathcal{H}^+$ and $\mathcal{H}^-$ are empty. We will show by contradiction
 that $\mathcal{H}^+=\emptyset,$ the case $\mathcal{H}^-=\emptyset$ is
 analogous.

 Fix a point $x_0\in \mathcal{H}^+$ which we consider as the basepoint
 of $M$ and put $\Gamma=\pi_1(M,x_0)$. We consider the action of
 $\Gamma$ on $Q^{\phi}$ which is induced by the action of
 $\Gamma$ on $\widetilde{M}$ by covering automorphisms.
 Let $\mathcal{G}^s$ and $\mathcal{G}^u$ be the foliations on
 $Q^{\phi}$ which are induced by $\widetilde{\cF}^s$ and $\widetilde{\cF}^u$,
 respectively. They are both preserved by the $\Gamma$-action.

 Let $\theta_0$ be the $\tilde{\phi}$-orbit of a lift of $x_0$ in $\widetilde{M}$.
 Since the $\phi$-orbit of $x_0$ is diffeomorphic to $\mathbb{T}^k$
 and incompressible, the isotropy group $\Gamma_0$ of $\theta_0$ is isomorphic
 to $\ZZ^k$. Let $F_0$ be the leaf through $\theta_0$ of $\mathcal{G}^s$
 and put $F_0'=F_0-\{\theta_0\}$.

 As the foliation $\mathcal{G}^u$ is orientable and one dimensional,
 their leaves admit a natural order. For all $x\in Q^{\phi}$ the subset
 of $\mathcal{G}^u(x)$ comprising elements above $x$ is
 denoted by $(x,+\infty)$.
 The fact that $x_0\in \mathcal{H}^+$ means:
 $$
 (\theta_0,+\infty)\cap \Gamma\cdot F_0=\emptyset.
 $$
 On the other hand, since $\theta_0$ is only point of $F_0$ which is the
 lift of a compact orbit, all the points of $F_0'$ are not lifts of
 compact orbits. This means that each point $x\in F_0'$ is not a
 lift of an orbit contained in $\mathcal{H}^+,$ equivalently,
 $(x,+\infty)\cap \Gamma\cdot F_0 \neq \emptyset.$ Let $h(x)$ be the
 infimum of $(x,+\infty)\cap \Gamma\cdot F_0$. We observe that $h$
 is, by definition, injective.

 \textbf{Claim 1.}
 \textit{Every $x\in F_0'$ is strictly inferior to $h(x)$.}
 Indeed, there exists $[x,\beta_x )$ a neighborhood of $x$
 in $[x,+\infty)$ such that all the leaves of $\mathcal{G}^s$
 which meets $[x,\beta_x )$ also meets $[\theta_0,+\infty )$.
 Hence, since none of these leaves is of the form $\gamma\cdot F_0,$
 $\gamma\in \Gamma,$ we obtain that
 $\beta_x \leq h(x)$. Consequently $h(x)>x.$

\begin{figure}[htb]
$$\vcenter{\hbox{\input{strsup.pstex_t}}} $$
\caption{\label{fig.beta}}
\end{figure}

 \textbf{Claim 2.}
 \textit{The image of $F_0'$ by $h$ is contained in a leaf $F_1$ of
  $\mathcal{G}^s$.} Given $x$, consider a small product neighborhood
  around $h(x)$. Then it is clear that $h(y)$ and $h(x)$ lies in the same
  leaf of $\cG^s$. Hence, for every weak stable leaf $F'$, the subset $\Omega(F')$
  of $F'_0$ comprising elements whose image by $h$ belongs to $F'$ is open.
  $F'_0$ is the disjoint union of all the $\Omega(F')$, and is connected (here we use the hypothesis
  $n > k+2$), therefore, all the $\Omega(F')$ are empty, except one, $\Omega(F_1)$.

 \textbf{Claim 3.}
 \textit{The leaf $F_1$ is $\Gamma_0$ invariant, and the map $h: F'_0 \to F_1$
 is a $\Gamma_0$-equivariant injective local homeomorphism onto its image.}
 It is clear that $h \circ \gamma = \gamma \circ h$ for every
 $\gamma$ in $\Gamma_0$. The claim follows.

 According to Claim 1, $F_1$ and $F_0$ are disjoint. Consider the stable leaf $\bar{F}_1=\pi(F_1)$:
 it is a bundle with contractible fibers over a flat cylinder $\RR^k/\Gamma_1$.
 Its fundamental group is $\Gamma_0$, hence $\Gamma_1$ is isomorphic to $\ZZ^k$, and thus,
 a lattice in $\RR^k$. According to Lemma~\ref{le:compactcausal}, it contains an Anosov element,
 which acts as a contraction in every strong stable leaf. Therefore, $\bar{F}_1$ contains
 a unique compact orbit. This compact orbit lifts as an element $\theta_1$, which is the unique
 $\Gamma_0$-fixed point in $F_1$. Let $F'_1 = F_1 - \{ \theta_1 \}$. Observe that since
 $h$ is injective, that $\theta_1$ is a $\Gamma_0$-fixed point,
 and that $\Gamma_0$ admits no fixed point in $F'_0$, the image of $h$ is contained in $F'_1$.

 \textbf{Claim 4.}
 \textit{The map $h: F'_0 \to F_1$ is a homeomorphism.} The only remaining point to show
 is the fact that $h(F'_0)=F'_1$.
 According to item (4) of Remark~\ref{rk.first},
 $\Gamma_0$ contains an element $\gamma_0$ such that, for some Anosov element $v$ of
 $\RR^k$ we have $\phi^v(x)=\gamma_0 x$ for every $x$ in $\theta_0$. It follows that
 the action of $\gamma_0$ on $F_0$ is contracting, admitting $\theta_0$ as its unique fixed point.
 Therefore, the action of $\gamma_0$ on $F'_0 \approx \RR^p - \{0\}$ is free, properly discontinuous,
 and the quotient space is diffeomorphic to $\SS^{p-1} \times \SS^1$. If we knew that
 the action of $\gamma_0$ on $F_1$ is also a contraction, then the claim would follow immediately
 from the fact that $h$ induces a continuous map between the quotient spaces $F'_0/<\gamma_0>$
 and $F'_1/<\gamma_0>$, and from the compactness of these quotient spaces.
 Unfortunately, there is no warranty that $\gamma_0$ acts properly on $F'_2$, hence we need
 a slightly more intricate argument.

 Since $\gamma_0$ is a contraction, there is
 an embedded codimension $1$ sphere $S_0$ in $F'_0$, boundary of a closed ball $B_0$ containing
 $\theta_0$, such that $\gamma_0(S_0)$ is another embedded sphere, contained inside $B_0$,
 and disjoint from $S_0$. The union $S_0 \cup \gamma_0(S_0)$ is the boundary of a subdomain
 $W_0 \approx \SS^{p-1}\times [0,1]$ of $B_0$. Now, the union of all the iterates
 $\gamma_0^nW_0$ covers the entire $F'_0$.

 Since $h$ is injective, its image is a $\Gamma_0$-invariant domain $W_\infty$ of
 $F'_1 \approx \SS^{p-1}\times\RR$ diffeomorphic to $\SS^{p-1}\times\RR$, containing
 the embedded sphere $S_1=h(S_0)$. Now observe that even if $\gamma_0$ might not be
 contracting in $F_1$, the same argument as the one used in $F'_0$ ensures that some element
 $\gamma_1$ of $\Gamma_0$ is contracting. Then, for $N$ sufficiently big, $\gamma_1^NS_1$
 is disjoint from $S_1$. By construction, $S_1$ does not bound a ball inside $W_\infty$, hence the
 same is true for $\gamma_1^NS_1$: there are both incompressible spheres inside
 $W_\infty \approx \SS^{p-1}\times\RR$. It follows that their union is the boundary of
 a compact domain $W_1 \subset W_\infty$. Considering $W_1$ as a compact domain in $F'_1$,
 we get that $F'_1$ is the union of the iterates under $\gamma_1^N$ of $W_1$. Therefore,
 $W_\infty$ is the entire $F'_1$.

\textbf{Conclusion.} Consider the sphere $S_0$ introduced in the previous step, and its image
$S_1$ by $h$. By construction, $S_1$ bounds a ball $B_0$ in $F_0$, containing $\theta_0$.
According to Jordan-Sch\"onflies Theorem, $S_1$ is also the boundary of a closed ball $B_1$ in
$F_1$. If $B_1$ does not contain $\theta_1$, then it would be contained in $F'_1$, and $h^{-1}(B_1)$
would be a closed ball in $F'_0$ bounded by $S_0$: contradiction.

Let $\cC$ be the union of all unstable segments $[x, h(x)]$ for $x$ describing $S_0$.
The union $\cS$ of $\cC$ with $B_0$ and $B_1$ is then a submanifold of $Q^\phi$, homeomorphic
to a sphere of codimension one. Since $Q^\phi$ is homeomorphic to $\RR^{n-k}$, $\cS$ is the boundary of
a closed topological ball $\cB$.

We not get the concluding final contradiction as follows: the $\cG^u$-leaf $\ell_1$
through $\theta_1$ is a closed
line in $Q^\phi$, crossing $\cS$ at $\theta_1$. Since $\cB$ is compact,
$\cG^u$ must escape from it, and thus, cross $\cS$ at another point. This intersection
cannot occur in $\cC$, since $\cC$ is tangent to $\cG^u$. It cannot occur in $B_1$,
since, as a leaf of $\cG^u$, it intersects every $\cG^s$-leaf in at most one point.
Therefore, $\ell_1$ must intersect $B_0$, and this intersection is reduced to one point.
Finally, since $F_0$ and $\ell_1$ are $\Gamma_0$-invariant, this intersection point must
be fixed by $\Gamma_0$: hence, it is $\theta_0$.

Therefore, $\ell_1$ contains two $\Gamma_0$-fixed points: $\theta_0$ and $\theta_1$.
This is a contradiction with the fact that unstable leaves contain at most one compact orbit
(see figure~\ref{fig.secoupe})

\begin{figure}[htb]
$$\vcenter{\hbox{\input{secoupe.pstex_t}}} $$
\caption{\label{fig.secoupe}}
\end{figure}

 \end{proof}

 \section{Conclusion}
 \label{sec:conclusion}
 As we already mentioned in the introduction, codimension one Anosov flows has been extensively
 studied, from the 60's until nowadays. It is reasonable to expect that all these results
 admit natural extensions to (irreducible) Anosov actions of $\RR^k$, but most work still has to be done.

 A \textit{symmetric flow} is the flow defined by a one-parameter subgroup $g^t$
 of a Lie group $G$ by right translations on a quotient manifold $\Gamma\backslash{G}/K$, where
 $\Gamma$ is a lattice of $G$ and $K$ a compact subgroup commuting with $g^t$.
 In \cite{tomter1} P. Tomter classified Anosov symmetric flows up to finite coverings and conjugacy
 when $G$ is semisimple or solvable. He proved that in the former case, the symmetric flow is
 (commensurable to) the geodesic flow of a rank $1$ symmetric space, and in the former case,
 the flow is (commensurable to) the suspension of hyperbolic automorphisms of
 a compact infranilmanifold. He further pursued his study to the more general case (\cite{tomter2}).

 This definition of symmetric flows extends naturally to the notion of symmetric actions of
 $\RR^k$. It is natural to ask about the classification of these actions for $k > 1$, at least
 in the case of irreducible actions. But the case $k=1$ is already quite intricate. In a forecoming
 paper, we will classify irreducible symmetric actions of $\RR^k$ of codimension one:
 either they are Anosov symmetric flows, or suspensions of hyperbolic automorphisms of tori (cf.
 examples~\ref{ex:suspension}, \ref{ex:diffeo}).

 In \cite{ghys}, E. Ghys proved that Anosov flows of codimension one on a manifold of dimension $\geq 4$,
 preserving a volume form
 and for which the sum of the stable and the unstable bundles is $C^1$, is topologically equivalent to
 the suspension of an Anosov diffeomorphism (and hence a hyperbolic automorphism of the torus).
 In a forecoming paper, we also extend this result to the $k \geq 2$ case.

 Actually, this last statement, for Anosov flows, has been recently highly improved: S. Simic
 proved that the same conclusion holds if the sum of the stable and unstable bundles is only
 Lipschitz regular (but this is still a restrictive hypothesis) (\cite{simic1}), and Asoako furthermore
 proved that the volume form preserving form hypothesis can be removed (\cite{asaoka}). All these
 impressive results are important steps towards the Verjovsky conjecture: \textit{every codimension one Anosov flow
 on a manifold of dimension $\geq 4$ is topologically equivalent to a suspension (of a hyperbolic toral
 automorphism).} Moreover, S. Simic announced a complete solution of Verjovsky's conjecture
 (\cite{simic2}).

 Therefore, it seems reasonable to conjecture:

 \vspace{.5cm}
 \noindent
 \textbf{Conjecture.}
 \textit{Every irreducible codimension one Anosov action of $\RR^k$ on a manifold of dimension $\geq k+3$
 is topologically conjugate the suspension of an Anosov action of $\ZZ^k$ on a closed manifold. }

%%%$\mathscr{A B C D E F G H I J K L M N O P Q R S T U V X Y W Z}$
%%%$\mathcal{F} \  F   $

\end{document}